\documentclass[review,3p]{elsarticle}

\usepackage{lineno}
\usepackage{algorithmic} 
\usepackage{algorithm}
\usepackage{float}
\usepackage{lipsum}
\usepackage{subfigure}
\usepackage{amssymb}
\usepackage{wrapfig}
\usepackage{graphicx}
\usepackage{color}
\usepackage{easylist}
\usepackage{amsmath}
\usepackage{easylist}

\renewcommand{\vec}[1]{\ensuremath\boldsymbol{#1}}

\newcommand{\gphi}{\nabla\phi}
\newcommand{\ngphi}{\left|\left|\nabla\phi\right|\right|}

\newcommand{\newtext}[1]{#1}

\newenvironment{lequation}[0]{\begin{linenomath}\begin{equation}}{\end{equation}\end{linenomath}}

\modulolinenumbers[1]

\renewcommand{\vec}[1]{\ensuremath\boldsymbol{#1}}

\journal{ Computer Methods in Applied Mechanics and Engineering}

\begin{document}

\begin{frontmatter}

\title{Level Set Jet Schemes for Stiff Advection Equations: The SemiJet Method}	

\author[UB]{Guhan Velmurugan}
\author[UNC]{Ebrahim M. Kolahdouz}
\author[UB]{David Salac\corref{cor1}}
\ead{davidsal@buffalo.edu}
\address[UB]{Mechanical and Aerospace Engineering, University at Buffalo SUNY, 318 Jarvis Hall, Buffalo, NY 14260-4400, USA}
\address[UNC]{Department of Mathematics, The University of North Carolina at Chapel Hill, Phillips Hall, Chapel Hill, NC 27599-3250, USA}
\cortext[cor1]{Corresponding Author}

\begin{abstract}
	Many interfacial phenomena in physical and biological systems are dominated by high order geometric quantities such as curvature. 
	Here a semi-implicit method is combined with a level set jet scheme to handle stiff nonlinear advection problems.
	The new method offers an improvement over the semi-implicit gradient augmented level set method previously introduced
	by requiring only one smoothing step when updating the level set jet function while still preserving the underlying methods higher
	accuracy. Sample results demonstrate that accuracy is not sacrificed while strict time step restrictions can be avoided.
\end{abstract}


\end{frontmatter}


\section{Introduction and Overview}
\label{sec:1.0}
Level set methods are one of the most popular interface capturing approaches 
with a wide variety of applications ranging from fluid flow \cite{Sethian2003}, image segmentation \cite{Malladi1995}, materials 
science \cite{chen1997simple} and structural topology optimization \cite{wang2003}.  
In the level set approach an interface $\Gamma(t)$ is
represented as a zero level set of a smooth function $\phi(\vec{x},t)$:
\begin{lequation}
	\Gamma(t)=\left\{\vec{x}:\phi(\vec{x},t)=0\right\}.
\end{lequation}
The evolution of the interface is implicitly tracked by evolving $\phi(\vec{x},t)$ due to an underlying velocity field $\vec{u}$,
\begin{lequation}
	\label{eq:level-set-advection}
	\frac{D\phi}{Dt}=\frac{\partial \phi}{\partial t}+\vec{u}\cdot\gphi=0.
\end{lequation}

Level set methods are robust in handling topological changes and also enable easy calculation of geometric quantities of
the interface. For instance, the unit normal pointing outward from the interface and the curvature of the interface can be expressed
in terms of $\phi(\vec{x},t)$:
\begin{linenomath}\begin{align}
	\label{eq:normal} \vec{n}=&\frac{\gphi}{\ngphi}, \\
	\label{eq:curvature} \kappa=&\nabla\cdot\frac{\gphi}{\ngphi}.
\end{align}\end{linenomath}

The accuracy of the level set method can be improved by including the gradients of the level set
and advecting them in a fully coupled manner. This method is known as the gradient augmented level set and 
was first introduced by Nave et al. as a globally third order accurate method 
that uses p-cubic Hermite interpolation with an optimally local stencil \cite{Nave2010}. Benefits of incorporating gradient information also 
include: 1) representation of subgrid structures, 2) simple and accurate approximation of the 
curvature and 3) simple implementation of adaptive mesh refinement with optimal locality \cite{kolomenskiy2014adaptive}. 
When used for multiphase flow these benefits result in increased volume conservation over the base level set method \cite{lee2014narrow}.

Jet schemes, developed by Seibold et al., are a natural extension of the original gradient augmented level set method and 
allow for the creation of schemes with an arbitrary high order of accuracy \cite{Seibold2011}. 
The higher order of accuracy is accomplished by tracking a suitable jet of the solution, i.e.
appropriate derivative information in addition to function values. The higher order spatial derivatives are utilized to construct
high order Hermite interpolations. Jet schemes of up to fifth order accuracy were developed and shown to obtain
comparable results to that of WENO5 \cite{Seibold2011}. Chidyagwai et al. demonstrated that jet schemes are computationally more efficient 
than comparable WENO schemes while possessing the optimal locality of discontinuous Galerkin methods \cite{Chidyagwai}.

Many problems concerning moving interfaces require computation of higher order spatial derivatives of the level set.
For example, the surface tension driven flows of immiscible fluids contain nonlinear and nonlocal higher order terms due to the Laplace-Young condition at the interface 
\cite{Sethian2003, Sussman1994, hysing2006new, salac2008local, popinet2009accurate}. 
As another example, in modeling the hydrodynamics of lipid vesicles the bending resistance exerted on the surrounding fluid 
is proportional to the second derivative of the curvature \cite{Salac2011}.  This involves the computation of up to the fourth order derivative of the 
level set function which makes the evolution of the interface front extremely stiff.
The gradient augmented level set method and jet schemes were originally developed for linear advection equations and are not able
to maintain a smooth level set field without strict time step restrictions \cite{kolahdouz2013semi}.

The semi-implicit formulation of an Eulerian level set method was first introduced by Smereka
and allows for time steps on the order of the grid spacing  \cite{Smereka2003}.
This method is based on 
extracting the linear portion of the advection
equation and treating that implicitly while using an explicit form of the nonlinear part.
Semi-implicit methods are a powerful alternative to explicit integration methods which suffer from 
strict time step conditions and fully implicit methods which require solving a non-linear system of equations for every time step \cite{salac2008local, Smereka2003}. 

A combination of the standard semi-implicit method and the original gradient augmented level set method was recently reported under the name
of the Semi-implicit Gradient Augmented Level
Set (SIGALS) \cite{kolahdouz2013semi}.
The method provides more accuracy through tracking the level set gradients and
alleviates some of the nonlocal behavior of the original semi-implicit level set method.
Another closely related formulation of this method 
has been successfully implemented in the simulation of vesicle electrohydrodynamics and the numerical experiments show
the robustness of the semi-implicit formulation in capturing the fast and long-range deformation of the vesicle in strong DC electric fields \cite{KolahdouzPRE2015, kolahdouz2014electrohydrodynamics}.
In another recent work a semi-implicit method was 
introduced to treat the numerical instabilities of partial differential equations in general \cite{duchemin2014explicit}.\par

In this work the level set jet scheme is extended to model interface evolution due to numerically stiff velocity fields. 
Unlike the SIGALS method, which requires four smoothing operations per time step for a three-dimensional problem and updates the gradients through an analytical differentiation method,
here a different approach is introduced to keep the level set and gradients smooth with less computational cost. Additionally, the derivatives are updated
through a simpler and more effective $\epsilon$-finite difference method. Compared to the original SIGALS method, 
which uses an analytic updating method for the derivative fields, this $\epsilon$-finite difference method extends more easily to higher 
level-set derivative fields. It has also been previously demonstrated that the accuracy of the $\epsilon$-finite difference method
and analytic method is identical, while the $\epsilon$-finite difference method has a lower computational cost \cite{Chidyagwai}.

The remainder of this paper is organized as follows. In Section \ref{sec:numback} the numerical background information is presented. This includes
the base level set jet scheme and the 
original level set semi-implicit formulation. In Section \ref{sec:semijet} the new semi-implicit jet scheme (SemiJet) is described.
Convergence and stability analysis using sample two and three dimensional results are provided in Section \ref{sec:results}.
A brief conclusion and future work is given in Section \ref{sec:conclusion}.

\section{Numerical Background}
\label{sec:numback}

The method presented in this paper is built upon the original level set jet scheme and a semi-implicit formulation for Hamilton-Jacobi equations.
In this section some background information on the original level set jet scheme and the semi-implicit level set method are given.
A brief mention of velocity extension and reinitialization is also presented.

\subsection{Level Set Jet Schemes}
\label{subsec:basejet}

Consider the advection of an interface using the level set method on a fixed grid. The level set function is only known at
grid points. To obtain interface information away from grid points some type of interpolation must be used. The idea
behind level set jet schemes is to advect some or all of the information needed to obtain high-order interpolation functions
in a local manner without the need to use wide stencils \cite{Seibold2011}. A grouping of this information, called a jet, consists of the 
level set function and certain derivatives of the level set. For example, using a jet which consists of the level set function, $\phi$,
in addition to gradient vector field, $\phi_x$ and $\phi_y$, and the first cross-derivative, $\phi_{xy}$, it would be possible to 
define a cubic Hermite interpolant on a two-dimensional Cartesian grid without the need for derivative approximations.

A critical component of any level set jet scheme is the coupled advection of the jet information. The original gradient-augmented
level set method and the more recent level set jet schemes both use a semi-Lagrangian approach: they evolve values 
on a fixed Eulerian grid by tracing characteristics using a Lagrangian method \cite{Nave2010, Seibold2011}. Consider first the 
advection of the level set function alone. In this case the time-evolution is given by Eq. (\ref{eq:level-set-advection}). Using the 
material derivative form of the evolution equation it can be written to first order that
\begin{linenomath}\begin{align}
	\frac{D\phi}{Dt}\approx \frac{\phi_{n+1}-\phi^d_n}{\Delta t}=0,
\end{align}\end{linenomath}
where $\phi_{n+1}$ is the level set value at a grid point $\vec{x}$ at time $t_{n+1}$, $\phi^d_n$ is the level set value at the departure (or foot)
location, and $\Delta t$ is the time step. The departure location can be obtained from
\begin{lequation}
	\frac{D\vec{x}(t)}{Dt}=\vec{u}(\vec{x}(t),t).
\end{lequation}
To first order this can be written as
\begin{lequation}
	\vec{x}^d=\vec{x}-\Delta t\;\vec{u}(\vec{x},t_n).
	\label{eq:depLocation}
\end{lequation}
The updated level set value can then be obtained by evaluating an interpolation function $P_{\phi}$,
\begin{lequation}
	\phi_{n+1}=\phi_n^d=P_{\phi}(\vec{x}^d, t_n).
	\label{eq:updatedValue}
\end{lequation}
Higher order method can be developed and will be discussed later.

A critical component of Jet schemes is the advection of the additional information. This must be done such that derivatives
of the level set do not de-couple from the base level set function. There are two basic methods to accomplish this. The first
is called analytic differentiation and was used in the gradient-augmented level set method \cite{Nave2010}. In the 
analytic differentiation method spatial derivatives of the level set evolution equation are used to 
evolve of the derivative fields. For example, to advect the gradient of the level set, $\vec{\psi}=\nabla\phi$, an advection equation 
can be obtained by taking the gradient of Eq. (\ref{eq:level-set-advection}),
\begin{lequation}
	\frac{D\vec{\phi}}{Dt}=-\nabla\vec{u}\cdot\vec{\psi}.
\end{lequation}

\begin{figure}[ht!]
	\begin{center}
		\includegraphics[width=0.25\textwidth]{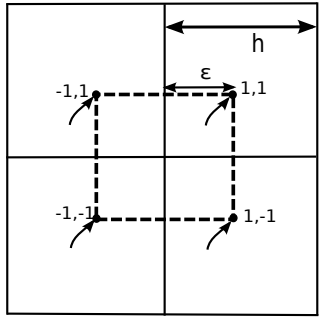}		
		\caption{The subgrid used in the $\epsilon$-finite difference update method for a Jet.}
		\label{fig:subgrid}
	\end{center}
\end{figure}

The analytic differentiation method does not work well with higher-order derivatives as the resulting evolution equations can become quite cumbersome.
An alternative is the $\epsilon$-finite difference method. At every grid point a local Cartesian sub-grid is constructed with a spacing of $\epsilon\ll h$, where
$h$ is the grid spacing. This defines a grid given by $\vec{x}^{\vec{q}}=\vec{x}+\vec{q}\epsilon$ for 
$\vec{q}\in\{-1,1\}^p$ where $p$ is the dimension of the simulation, either $p=2$ or $p=3$, Fig. \ref{fig:subgrid}.
At each subgrid point an updated level set value, $\phi^{\vec{q}}$, is obtained using Eqs. (\ref{eq:depLocation}) and (\ref{eq:updatedValue}).

The updated level set values on the subgrid can be used to update the derivatives at grid points by using finite difference approximations.
In two dimensions values up to $\phi_{xy}$ can be obtained via the following finite-difference approximations:
\begin{linenomath}\begin{align}
	\label{eq:jetPhi} \phi = &\frac{1}{4}\left(\phi^{(1,1)}+\phi^{(-1,1)}+\phi^{(1,-1)}+\phi^{(-1,-1)}\right), \\
	\label{eq:jetPhix} \phi_x = &\frac{1}{4\epsilon}\left(\phi^{(1,1)}-\phi^{(-1,1)}+\phi^{(1,-1)}-\phi^{(-1,-1)}\right), \\	
	\label{eq:jetPhiy} \phi_y = &\frac{1}{4\epsilon}\left(\phi^{(1,1)}+\phi^{(-1,1)}-\phi^{(1,-1)}-\phi^{(-1,-1)}\right), \\
	\label{eq:jetPhixy} \phi_{xy} = &\frac{1}{4\epsilon^2}\left(\phi^{(1,1)}-\phi^{(-1,1)}-\phi^{(1,-1)}+\phi^{(-1,-1)}\right).
\end{align}\end{linenomath}

The overall accuracy of the scheme is $\mathcal{O}(\epsilon^2)+\mathcal{O}(\delta)$ for Eq. (\ref{eq:jetPhi}), 
$\mathcal{O}(\epsilon^2)+\mathcal{O}(\delta/\epsilon)$ for Eqs. (\ref{eq:jetPhix}) and (\ref{eq:jetPhiy}), 
and $\mathcal{O}(\epsilon^2)+\mathcal{O}(\delta/\epsilon^2)$ for Eq. (\ref{eq:jetPhixy}), where $\delta$ is machine precision \cite{Seibold2011}.
An optimal choice is thus $\epsilon=\mathcal{O}(\delta^{1/4})$, which yields an overall accuracy of $\mathcal{O}(\delta^{1/2})$ in all cases.
Note that this particular choice for $\epsilon$ ensures that the overall errors are dominated by other terms, such as the time discretization,
and not the calculation of derivatives fields during updating.
While, in general, the derivatives in a level-set description of an interface are discontinuous in a finite number of points for 
a discretized domain, these discontinuities are typically away from the interface and situations where they occur near the 
interface are temporary for moving interfaces. These discontinuities have not been observed to be detrimental to the overall
accuracy of Jet schemes in the past \cite{Seibold2011,Chidyagwai}.

The major advantage of the $\epsilon$-finite difference method over the analytic differentiation method is the applicability to higher-order derivatives.
Higher-order derivative values can be calculated by expanding the subgrid stencil and using appropriate finite difference approximations.
The small subgrid size ensures that results are accurate. Additionally, Chidyagwai et al.
showed that $\epsilon$-finite difference method are faster than analytical differentiation and that 
the method is more efficient than WENO while possessing the accuracy of discontinuous Galerkin schemes \cite{Chidyagwai}.

\subsection{The Semi-Implicit Level Set Method}
\label{subsec:basesemi}

One issue with all interface tracking techniques is their instability if the underlying velocity field is numerically stiff.
A classic example is mean curvature flow where the velocity of the interface is given by 
$\vec{u}=-\kappa\vec{n}$. From Eq. (\ref{eq:curvature}) it is clear that the curvature is a second-order derivative
of the level set function, which would require a time step restriction of $\Delta t=\mathcal{O}(h^2)$ for numerical stability if
an explicit time discretization is used.

The semi-implicit level set method was developed to circumvent this time step restriction \cite{Smereka2003}. 
If the evolution is relatively simple, such as mean-curvature flow, then it is possible to extract the linear portion 
of the level set evolution equation. This linear portion would be done implicitly while the remaining non-linear portion is done explicitly:
\begin{lequation}    
	\frac{\phi_{n+1}-\phi_{n}}{\Delta t}=\nabla^2\phi_{n+1}-\frac{\gphi_{n}}{\|\gphi_{n}\|}\cdot\nabla(\|\gphi_{n}\|).
\end{lequation}
In cases where it is not feasible to explicitly extract the linear portion a smoothing operator can be added to the 
evolution equation \cite{salac2008local, Smereka2003, duchemin2014explicit},
\begin{lequation}
	\frac{\phi_{n+1}-\phi_n}{\Delta t}+\vec{u}_n\cdot\nabla\phi_n-L(\phi_{n+1})+L(\phi_n)=0,
	\label{eq:SemiImplicit}
\end{lequation}
where the $L(\phi)$ is a linear damping operator which can be chosen based on the stiffness of the velocity field $\vec{u}$.
If the order of the velocity field is even, a common choice for the damping operator is $L(\phi)=\beta\nabla^m\phi$, 
where $m$ matches the order of the velocity field and $\beta$ is a constant. 
Note that this stabilization is added after normalization and therefore $\beta$ is dimensionless. If the system is dimensional, then this $\beta$ is akin
to a diffusivity parameter. 
If the velocity order is odd, then a more appropriate choice 
is the Hilbert transform of the $m^{th}$-order spatial derivative, $L(\phi)=\beta\mathcal{H}(\nabla^m\phi)$.
This transform ensures that the right scaling is 
observed in the frequency domain~\cite{duchemin2014explicit}.

An explanation of why Eq. (\ref{eq:SemiImplicit}) aids in stability has been provided in Ref. \cite{duchemin2014explicit} and is briefly explained here. 
Consider the simpler problem of $u_t=\alpha u_{xx}$ on a one-dimensional domain with periodic boundary conditions.
As the system is of second-order, an appropriate linear damping operator has the form of $L(\phi)=\beta u_{xx}$, so the modified evolution equation 
is now $u_t-\beta u_{xx}+\beta u_{xx}=\alpha u_{xx}$. Consider the perturbation of a single Fourier mode,
$u(x_j, t_n)=\xi^{n \Delta t} e^{i k j h}$, where $\xi$ is the amplification factor, $h$ is the grid spacing, $\Delta t$ is the 
time step, and $k$ is the Fourier mode.
Using the perturbation in a first-order in time and second-order in space discretization and simplifying results in the 
following expression for the amplification factor,
\begin{lequation}	
	\xi^{\Delta t}=1 + \dfrac{\tfrac{\alpha \Delta t}{h^2}\left(e^{i k h}+e^{-i k h} -2\right)}{1-\tfrac{\beta\Delta t}{h^2}\left(e^{i k h}+e^{-i k h} -2\right)}.
\end{lequation}
Using the relation $e^{i k h}+e^{-i k h} -2=-4\sin^2(h k/2)$, the amplification factor is now
\begin{lequation}	
	\xi^{\Delta t}=1 - \dfrac{\tfrac{4\alpha \Delta t}{h^2}\sin^2(h k/2)}{1+\tfrac{4\beta\Delta t}{h^2}\sin^2(h k/2)}.
\end{lequation}
A stable scheme is one where $|\xi^{\Delta t}|\leq 1$, which corresponds to the relation
\begin{lequation}	
	\left(4 \alpha-8\beta\right)\dfrac{\Delta t}{h ^2} \sin^2\left(\dfrac{k h}{2}\right)\leq 2.
\end{lequation}
As the $\sin^2$ term has a maximum value of 1, stability is ensured when
\begin{lequation}
	\left(4 \alpha-8\beta\right)\dfrac{\Delta t}{h ^2}\leq 2.
\end{lequation}
When $\beta=0$, which corresponds to no stabilization, the standard time-step restriction for second-order differential 
equations is obtained, $\alpha \Delta t/h^2\leq 1/2$. Any value of $\beta>0$ stabilizes the scheme as the condition now 
becomes $\alpha \Delta t/h^2\leq 1/2+2 \beta \Delta t/h^2$. In this specific case, when $\beta=\alpha/2$, the system becomes 
is equivalent to a Crank-Nicolson discretization and is unconditionally stable. 
For general advection, the value of $\beta=1/2$ has been shown to work 
well in most cases~\cite{Smereka2003}. Additionally, it has been shown that any scheme can be 
stabilized, even if the linear damping operator is of lower order than the velocity field~\cite{duchemin2014explicit}.

Finally, a word should be said about the boundary conditions for the semi-implicit update. If the domain is periodic then no additional
boundary conditions are required. For non-periodic boundaries a logical choice is to use homogeneous Neumann boundary conditions \cite{Burger2005,Olsson2005}. 
In the case of complex or curved boundaries it would be necessary to use an interface capturing technique \cite{Kublik2013,Liu2000}. 
Note that if level set reinitialization is used, as described below, the particular choice of boundary conditions can play little role 
in the evolution, as the level set function far away from the interface is replaced with a different, but equivalent one.

\subsection{Level Set Reinitialization and Extension}
\label{subsec:extension_reinit}

While any level set function where a contour describes the interface can be chosen, a typical choice is to use a signed-distance function
where the interior is given by $\phi<0$ while the exterior has $\phi>0$ and $|\phi(\vec{x})|$ indicates the distance from 
$\vec{x}$ to the interface. Even if a signed-distance function is chosen as the initial level set it can not be
guaranteed that it will remain so over time. Typically level set reinitialization is performed periodically 
to return the level set back to a signed distance function. One possible method is to use a PDE-based reinitialization
technique \cite{Sussman1994,Peng1999}, 
\begin{lequation}
	\frac{\partial \phi}{\partial \tau}+\textrm{sgn}(\phi_0)\left(\|\nabla\phi\|-1\right)=0,
\end{lequation}
where $\tau$ is a pseudo-time and $\textrm{sgn}(\phi_0)$ is the sign function
of the original level set. Another possible method is to reinitialize by explicitly calculating the closest point
from a point in the vicinity of the interface to the interface \cite{chopp2001}. While slower computationally, this explicit calculation
of the closest point has the advantage of higher accuracy than the PDE-based method and will have additional 
benefits for Jet schemes, as will be explained later.

In many situations the velocity is only defined at the interface. To advance a level set function
the velocity needs to be determined, at a minimum, in a region surrounding the interface. 
An extension algorithm can accomplish this by extending a function in the direction normal to the interface. One method to do this 
is to solve a hyperbolic PDE in the form of \cite{Peng1999}
\begin{lequation}
	\frac{\partial q}{\partial \tau}+\textrm{sgn}(\phi)\vec{n}\cdot\nabla q=0,
\end{lequation}
where $q$ is the quantity to be extended from the interface. This PDE can be solved for a few iterations using a pseudo-time step of $\Delta \tau\approx h/2$ to extend
the quantity $q$ in a region near the interface. 
An alternative to the PDE-based extension algorithm is to explicitly extend a quantity by interpolating 
the quantity $q$ at the closest point which could have been determined during a reinitialization procedure. In practice both methods work well.

\section{The Semi-Implicit Level Set Jet Scheme (SemiJet)}
\label{sec:semijet}

\begin{figure}[ht!]
	\begin{center}
		\subfigure[]{
			\includegraphics[width=0.25\textwidth]{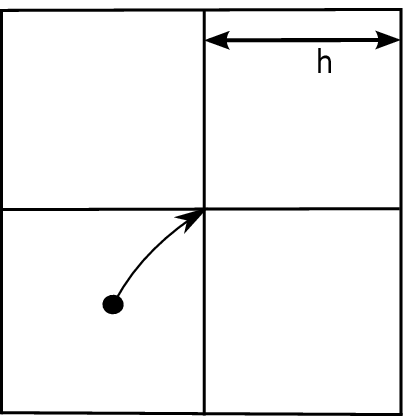}
			\label{fig:grid}	
		}
		\subfigure[]{
			\includegraphics[width=0.25\textwidth]{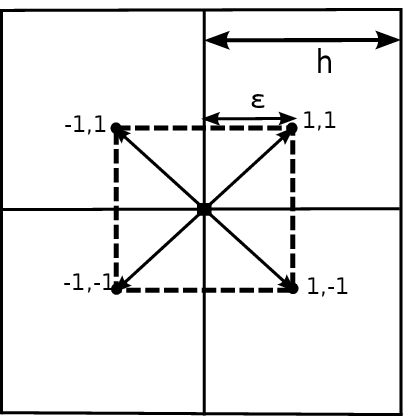}
			\label{fig:grid2}	
		}
		\subfigure[]{
			\includegraphics[width=0.25\textwidth]{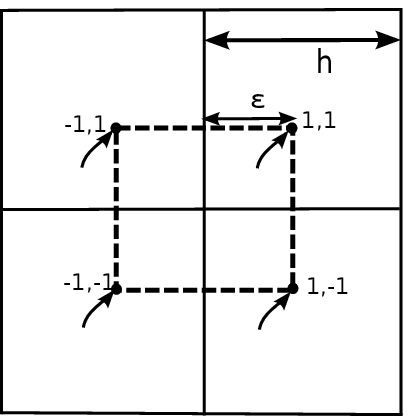}
			\label{fig:grid3}	
		}	
		\caption{The three basic steps of the SemiJet method. (a) A standard semi-Lagrangian semi-implicit step which updates grid points. (b) The effect of smoothing
		is applied to the sub-grid locations by way of a smoothing source term. (c) The subgrid points are updated using a semi-Lagrangian advection with the smoothing 
		source term, allowing for the remaining Jet fields to be updated.}
		\label{fig:grids}		
	\end{center}
\end{figure}

The idea behind the SemiJet method is to combine the increased accuracy of the standard level set jet method with the stability of the semi-implicit
level set method. The general outline has three steps: 1) a standard semi-Lagrangian, semi-implicit update, 2) determining the effect of smoothing, and 
3) propagation of the smoothing effect to derivative updates, see Fig. \ref{fig:grids} for a graphical representation. 
\newtext{Without loss of generality, only even-ordered velocity field will be considered. Therefore, the linear damping operator will have the form of 
$L(\phi)=\beta \nabla^m \phi$. If another linear damping operator is needed, such as when the velocity field is dominated by odd-order derivatives,
then $\nabla^m\phi$ is simply replaced by an appropriate operator.}

The first step is a standard semi-Lagrangian semi-implicit step on the base level set function $\phi$. Begin by writing the semi-implicit
level set equation in semi-Lagrangian form,
\begin{lequation}
	\frac{D\phi}{Dt}-\beta \nabla^m\phi_{n+1}+\beta \nabla^m\hat{\phi}=0,
	\label{eq:lagrangian_semi} 
\end{lequation}
where $\hat{\phi}$ is an approximation to $\phi_{n+1}$.
A first-order discretization in time results in
\begin{lequation}
	\frac{\phi_{n+1}-\phi^d_n}{\Delta t}-\beta \nabla^m\phi_{n+1}+\beta \nabla^m\phi_n=0,
	\label{eq:phi_first_order}
\end{lequation}
with $\hat{\phi}=\phi_n$.
This can be solved through a multi-step process given by
\begin{linenomath}\begin{align}   
	\label{eq:departure_location} \vec{x}^d &= \vec{x}-\Delta t\;\vec{u}(\vec{x},t_n), \\
	\label{eq:departure_phi_first} \phi_n^d &= P_{\phi}\left(\vec{x}^d,t_n\right), \\
	\label{eq:smoothing_first_order} \frac{\phi^{n+1}-\phi_n^d}{\Delta t}&=\beta \nabla^m \phi^{n+1} - \beta \nabla^m \phi^n,
\end{align}\end{linenomath}
where $\phi_n^d$ is the tentative (departure) level set value and $P_\phi\left(\vec{x}^d,t_n\right)$ is the interpolant
of the level set value at time $t_n$ evaluated at the departure location $\vec{x}^d$.

A second-order method can be obtained by writing
\begin{linenomath}\begin{align}	
	\frac{3\phi_{n+1}-4\phi^d_n+\phi^d_{n-1}}{2\Delta t} - \beta \nabla^m \phi_{n+1} + \beta \nabla^m \hat{\phi}=0,
	 \label{eq:smoothing_second_order}
\end{align}\end{linenomath}
where $\hat{\phi} = 2\phi_n-\phi_{n-1}$.
The departure value of the level set at time $t_n$ is calculated by
\begin{linenomath}\begin{align}
    \hat{\vec{u}}_{n+1} &=2 \vec{u}_n-\vec{u}_{n-1}, \\
	\vec{x}_1 &= \vec{x}-\Delta t\;\hat{\vec{u}}_{n+1}, \\
    \label{eq:second_order_departure} \vec{x}^d_n &= \vec{x}-\frac{1}{2}{\Delta t}\left(\hat{\vec{u}}_{n+1}+ \vec{P}_{\vec{u}}\left(\vec{x}_1,t_n\right)\right), \\
    \phi^d_n&=P_{\phi}\left(\vec{x}^d_n, t_n\right),
\end{align}\end{linenomath}
where $\vec{P}_{\vec{u}}$ is the interpolation function of the velocity field. The departure value at time $t_{n-1}$ can then
be obtained with
\begin{linenomath}\begin{align}
    \vec{x}^d_{n-1} &= \vec{x}-2 {\Delta t} \vec{P}_{\vec{u}}\left(\vec{x}^d_n, t_n\right), \\
    \phi^d_{n-1} &= P_{\phi}\left(\vec{x}^d_{n-1}, t_{n-1}\right).
\end{align}\end{linenomath}
A smooth level set field can now be obtained by solving Eq. (\ref{eq:smoothing_second_order}). Note that, while it is possible to use higher order
time discretizations, the results presented below use second-order accurate approximations to the smoothing field, $\nabla^m$. In many cases the velocity 
field is normalized to be $\mathcal{O}(1)$, and it is typically not advisable to have the interface move more than one grid spacing per time step.
Therefore, a situation where $\Delta t\approx h$, where $h$ is the grid spacing is a standard situation and thus it is not advantageous to use a higher-order
time discretization as the method would be limited by the spatial discretization.

At this point a smooth and updated level set field has been determined but the remaining derivatives fields in the Jet have not yet been updated.
The application of the smoothing step with $\beta>0$ has the effect of modifying the evolution of $\phi$, with the magnitude of this deviation from 
the $\beta=0$ (e.g. no smoothing) case given by $\beta \nabla^m\phi_{n+1} - \beta \nabla^m\hat{\phi}$. Conceptually,
if this deviation was known at the beginning of a time step, there would be no need to perform a semi-implicit update, as the corrections 
needed to maintain a smooth interface would already be known.

This concept is used to update the remaining 
fields in the Jet. Capture the effect of smoothing
by defining
a smoothing source term,
\begin{lequation}
	S_{\phi}=\beta \nabla^m\phi_{n+1}-\beta \nabla^m\hat{\phi},
	\label{eq:smoothing_source}
\end{lequation}
where both $\phi_{n+1}$ and $\hat{\phi}$ have already been determined. 
The method now proceeds in a similar manner to the standard level set Jet scheme, with the only modification that the 
advection now has a source term. At each grid point construct a sub-grid of spacing $\epsilon\ll h$:
$\vec{x}^{\vec{q}}=\vec{x}+\vec{q}\epsilon$ for $\vec{q}\in\{-1,1\}^p$, see Fig. \ref{fig:grid2}. 
At each sub-grid location calculate an updated level set value, $\phi^{\vec{q}}$. For example, a first-order scheme 
would use the following steps:
\begin{linenomath}\begin{align}
	\vec{x}^d &= \vec{x}^{\vec{q}}-\Delta t\;\vec{P}_{\vec{u}}(\vec{x}^{\vec{q}},t_n), \\
	\phi_n^d &= P_{\phi}\left(\vec{x}^d,t_n\right), \\
	\frac{\phi^{\vec{q}}-\phi_n^d}{\Delta t}&=P_{S}\left(\vec{x}^{\vec{q}}\right),
\end{align}\end{linenomath}
where $P_{S}$ is the interpolant of the smoothing source term given in Eq. (\ref{eq:smoothing_source}). The second-order
scheme proceeds similarly,
\begin{linenomath}\begin{align}	
	\hat{\vec{u}}_{n+1} &= 2\vec{P}_{\vec{u}}\left(\vec{x}^{\vec{q}},t_n\right)-\vec{P}_{\vec{u}}\left(\vec{x}^{\vec{q}},t_{n-1}\right),	\\
	\vec{x}_1 &= \vec{x}^{\vec{q}}-\Delta t\;\hat{\vec{u}}_{n+1}, \\
    \vec{x}^d_n &= \vec{x}^{\vec{q}}-\frac{1}{2}{\Delta t}\left(\hat{\vec{u}}_{n+1}+ \vec{P}_{\vec{u}}\left(\vec{x}_1,t_n\right)\right), \\
	\vec{x}^d_{n-1} &= \vec{x}^{\vec{q}}-2 {\Delta t} \vec{P}_{\vec{u}}\left(\vec{x}^d_n, t_n\right), \\
    \phi^d_{n-1} &= P_{\phi}\left(\vec{x}^d_{n-1}, t_{n-1}\right),\\
	\frac{3\phi^{\vec{q}}-4\phi^d_n+\phi^d_{n-1}}{2\Delta t} &= P_{S}\left(\vec{x}^{\vec{q}}\right).	 
\end{align}\end{linenomath}
In both the first- and second-order cases the smooth source term is evaluated at the sub-grid point $\vec{x}^{\vec{q}}$.
Once the sub-grid level set values have been determined the updated derivative fields at grid points 
can be calculated using finite difference approximations such as Eqs. (\ref{eq:jetPhix})-(\ref{eq:jetPhixy}), or similar.
Note that in the SemiJet method the level set function is not given by Eq. (\ref{eq:jetPhi}), but is instead the result 
of the smoothing step, either Eq. (\ref{eq:smoothing_first_order}) or (\ref{eq:smoothing_second_order}).

A summary of the SemiJet method can be written as:
\begin{enumerate}
	\item Determine a tentative updated level set field using a standard semi-Lagrangian level set method.
	\item Solve the linear system arising from a semi-implicit update. This produces the updated level set field.
	\item Calculate the effect of smoothing by way of a smoothing source term.
	\item Use the smoothing source term to determine updated level set values on the sub-grid.
	\item Update the remaining fields in the Jet using finite difference approximations on the sub-grid.	
\end{enumerate}
From this summary is it clear that only one linear system needs to be solved per time step. This is an 
advantage over the original SIGALS method which required a linear system for each component of the Jet.

\section{Numerical Results}
\label{sec:results}
In this section sample two- and three-dimensional results using the method described in Sec. \ref{sec:semijet} are presented.
Unless otherwise specified a constant of $\beta=0.5$ is chosen and second-order accurate time integration is used.
In this work only curvature flows are considered and thus the smoothing operator is set to the Laplacian, $m=2$.
The sub-grid spacing was set to $\epsilon=10^{-4}$ while periodic boundary conditions were used in all cases.
The smoothing operator is discretized using second-order accurate isotropic finite differences \cite{Kumar2004,Patra2006}.
The linear system of equations resulting during updating
the level set function are solved using Flexible Generalized Minimal Residual (FGMRES) method with Geometric Algebraic Multi-Grid (GAMG) as the preconditioner,
both provided by the PETSc library \cite{petsc-web-page,petsc-user-ref,petsc-efficient}, with a relative tolerance of the solver set to $10^{-8}$.

All required interpolation functions are cubic interpolants. The interpolation 
coefficients for the velocity field, $\vec{P}_{\vec{u}}$, and smoothing source term, $P_S$, are determined using a standard 16-node
stencil in two-dimensions and 64-node stencil in three-dimensions. The stencil for the interpolation coefficients for the Jet, $P_{\phi}$, are determined based on
the order of the Jet. 
The Jet structures here will be referred to as either a 0-Jet or Partial 1-Jet (P1-Jet). These names refer to the number of derivative fields 
which are available for interpolation functions. A 0-Jet consists of just the level set field, $\phi$. Thus a 0-Jet will require 
either 16- or 64-nodes, depending on the dimension of the simulation, to determine
the cubic interpolant. A P1-Jet includes the gradient field, $\vec{\psi}= (\phi_x, \phi_y, \phi_z)$.
To compute the Jet interpolation function the first cross-derivative are also required. For a P1-Jet these are computed using 
the cell-based approximation \cite{Nave2010}. In two-dimensions the values of $\phi_{xy}$ are first calculated at edge-centers. These are then
interpolated to the cell corners (the grid points) to obtain the necessary value for the interpolant. 
In this way only information from the local cell is required
to determine the interpolation function of a P1-Jet. Additional information about this interpolation procedure can be found in the original 
gradient-augmented level set method \cite{Nave2010}.
An extension of this nomenclature would be a Full 1-Jet (1-Jet), which
would include both the gradient field and the first cross-derivatives. Both the P1-Jet and 1-Jet give similar results and 
therefore only the P1-Jet is included. Please note that this naming convention differs from Seibold et al. \cite{Seibold2011}.

Reinitialization is performed during every time step by explicitly calculating the closest point in a band around the interface.
This band consists of all nodes within four grid spacings of the interface.
A reinitialized Jet is given by
\begin{linenomath}\begin{align}
	\phi(\vec{x})&=\textrm{sgn}(\phi_0(\vec{x}))\|\vec{x}-\vec{x}_{\Gamma}\|, \\
	\vec{\psi}(\vec{x})&=\textrm{sgn}(\phi_0(\vec{x}))\frac{\vec{x}-\vec{x}_{\Gamma}}{\|\vec{x}-\vec{x}_{\Gamma}\|},
\end{align}\end{linenomath}
where $\phi_0$ is the original level set function and $\vec{x}_{\Gamma}$ is the closest point on the interface. The rationale 
for this particular reinitialization scheme over a PDE based one is the direct calculation of the updated gradient field. 
Additionally, it has been shown than the direct reinitialization method is fourth-order accurate for the level 
set field and third-order accurate for the gradient field if the Jet field is smooth \cite{Anumolu2013}. 
If discontinuities in the gradient field exist, the accuracy of the level set is locally reduced to second-order while that of the gradient to first-order.
In a moving interface it is expected that these type of discontinuities are temporally temporary.
Therefore, this reinitialization should not be detrimental to the overall scheme, 
which is expected to be second-order accurate in space. \newtext{The use of other schemes, such as that of Cheng and Tsai~\cite{Cheng2008},
could be extended to reinitialize jet methods, but that is beyond the scope of this work.}

Using Jet schemes there are two choices when calculating curvatures. The first is the direct calculation of the curvature at closest points
using derivatives of the cubic interpolant. The second is to calculate curvatures at grid points and then interpolate to the interface. 
During the course of the numerical experiments the stability of the first method was much worse than the second, and thus the first method was not pursued.
Therefore, curvatures are first calculated at the grid points by writing Eq. (\ref{eq:curvature}) as
\begin{linenomath}\begin{align}
	\kappa = \frac{\phi_{xx}\phi_y^2+\phi_{yy}\phi_x^2-2\phi_{xy}\phi_x\phi_y}{\left(\phi_x^2+\phi_y^2\right)^{3/2}},	
\end{align}\end{linenomath}
in two-dimensions with a similar expression for three-dimensions. For the 0-Jet all derivatives are computed using isotropic finite differences.
For results using a P1-Jet the best results were obtained by averaging derivatives. For example, when computing $\phi_x$, $\phi_{xx}$ and $\phi_{xy}$ 
for use in the curvature calculation the following averaging would be used:
\begin{linenomath}\begin{align}
	\phi_x &= \tfrac{1}{2}\left(\psi^x + D_x\phi\right),\\
	\phi_{xx} &= \tfrac{1}{2}\left(D_x\psi^x + D_{xx}\phi\right),\\
	\phi_{xy} &= \tfrac{1}{3}\left(D_y\psi^x + D_x\psi^y + D_{xy}\phi\right),
\end{align}\end{linenomath}
where $D_x$, $D_y$, $D_{xx}$ and $D_{xy}$ are the isotropic finite difference approximations for the derivatives and $\psi^x$ and $\psi^y$ are the components of the 
gradient field tracked by the Jet. While this does remove some of the locality of the Jet method it was found that the accuracy and stability were increased. This
also adds additional coupling between the components of the Jet.

\subsection{Mean Curvature Flow}
\label{subsec:mean-curvature-flow}
Flow by mean curvature is a standard test problem for assessing the stability and accuracy of interfaces moving under a numerically 
stiff velocity field. In this case the velocity field is given by $\vec{u}=-\kappa\vec{n}$, where $\kappa$ is 
the mean curvature and $\vec{n}$ is the outward-facing unit normal vector. 

\begin{figure}[ht!]
	\begin{center}
		\subfigure[]{
			\includegraphics[width=0.4\textwidth]{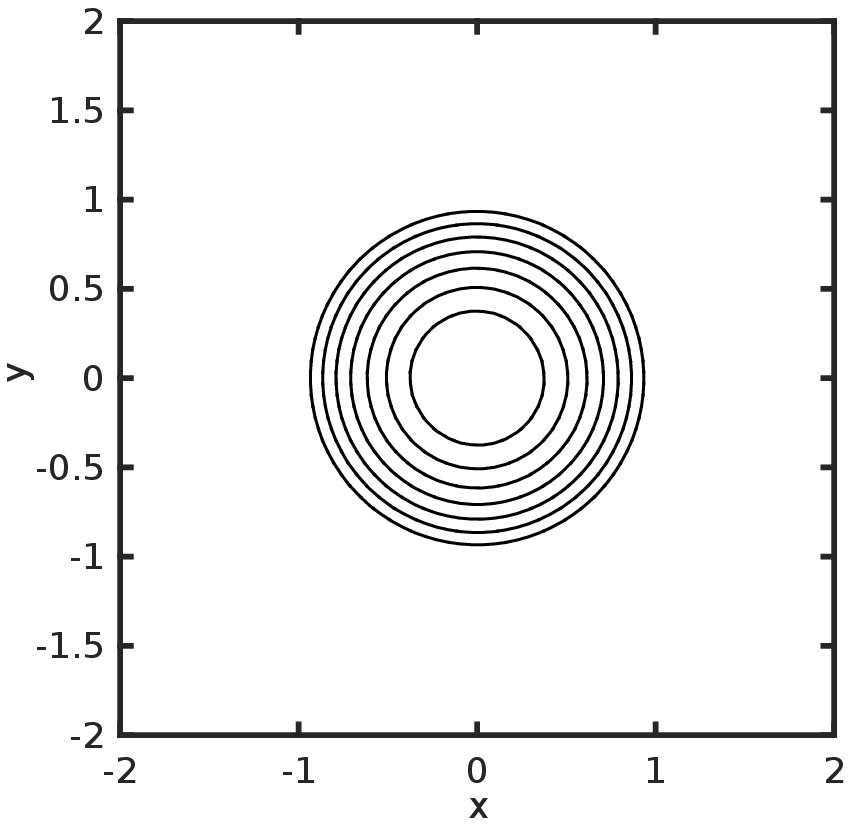}
			\label{fig:circleK}	
		}
		\subfigure[]{
			\includegraphics[width=0.4\textwidth]{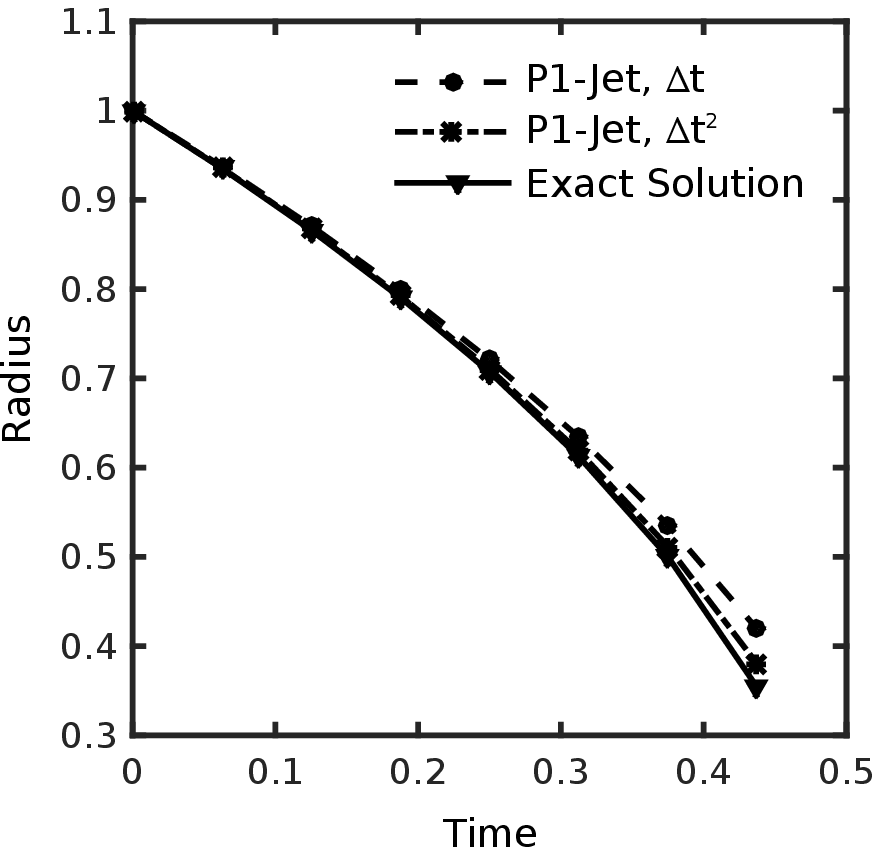}
			\label{fig:radiusvserror}	
		}	
		\caption{A circle collapsing under mean curvature flow using a P1-Jet with $h=0.125$ and $dt=4h^2$. (a) The interface location is shown at 
				intervals of $\Delta t$ until $t=0.375$.
			(b) A comparison between the radius over time using a P1-Jet with different time discretization methods.}	
		\label{fig:circle_plot}		
	\end{center}
\end{figure}

First consider the benchmark problem of curvature flow applied to a two dimensional circle of radius $r_0=1$ centered at $(0,0)$ in a $[-2,2]^2$ domain. 
The initial level set field is given by $\phi(x,y,0)=\sqrt{x^2+y^2}-r_0$, which is a signed distance function.
Under mean curvature flow the interface will shrink uniformly where the 
radius at any point in time is given by $r(t)=\sqrt{r_0^2-2t}$. 
An example of this behavior is provided in Fig. \ref{fig:circleK} using a grid spacing of $h=0.125$ and a time step of $\Delta t=4h^2$. 
A comparison of the average radius as calculated using the closest points
using a first-order and second-order in time scheme is presented in Fig. \ref{fig:radiusvserror}. The second-order method produces 
more accurate results, as should be expected.

\subsubsection{Convergence Study}
\label{subsubsec:convergence}
The overall accuracy of the method is a combination of the semi-Lagrangian and semi-implicit methods,
\begin{lequation}
	\label{eq:error}
	\mathcal{O}(\Delta t^k)+\mathcal{O}\big(\frac{h^{p+1}}{\Delta t}\big)+\mathcal{O}(h^L),
\end{lequation}
where $\Delta t$ is the time step, $h$ is the grid spacing, $k$ is the time step order, $p$ is the order of the interpolant, and $L$ 
is the accuracy order at which the smoothing operator is discretized. Here $k=1$ or $k=2$, $p=3$, and $L=2$. Clearly, there will be 
trade-off in accuracy as the time step is reduced. Given a fixed grid spacing it should be expected that at large time steps the error will be 
dominated by the time-integration error, $\mathcal{O}(\Delta t^k)$, and as such should decrease at a rate dictated by the time step order. 
As the time step decreases the error begins to be dominated by the smoothing term, $\mathcal{O}(h^L)$, if $L<p+1$.
A further decrease in the time step results in the interpolation error, $\mathcal{O}(h^{p+1}/\Delta t)$, dominating
and any further decrease in the time step will result in an increase in the error.

Therefore there are three convergence regions with respect to the time step which are to be expected. The first is a decrease 
at the rate of $\mathcal{O}(\Delta t^k)$, the second will be a relatively constant error region dictated by $\mathcal{O}(h^L)$, and the 
final will be an increase in the error at a rate of $\mathcal{O}(h^{p+1}/\Delta t)$.
Note that the $1/\Delta t$ is included to account for the fact that interpolation is performed about $1/\Delta t$ times during the course of a simulation.
This is not a unique occurrence and has been observed in other instances of semi-Lagrangian time integration \cite{strain1999semi, xiu2001semi}.

To verify this accuracy a two-dimensional circle shrinking under mean curvature flow is compared to the analytic solution. A circle with an
initial radius of 1 centered at the origin is allowed to evolve until a time of $T=0.375$, at which point the error is calculated.
Define the error at a point on the interface as $e(x_{\Gamma}, y_{\Gamma}, t)=|\sqrt{x_{\Gamma}^2+y_{\Gamma}^2}-r(t)|$, where 
$(x_{\Gamma}, y_{\Gamma})$ is a point on the interface and $r(t)=\sqrt{r_0^2-2t}$ is the analytic radius. The overall error is 
the  $L_\infty$-error defined over all closest points calculated at $T=0.375$. This error definition tests not only the 
time-evolution of the Jet but also the recovery of the interface by way of the cubic interpolation function.

\begin{figure}[ht!]
	\begin{center}
		\subfigure[$h=0.0625$]{
			\includegraphics[width=0.4\textwidth]{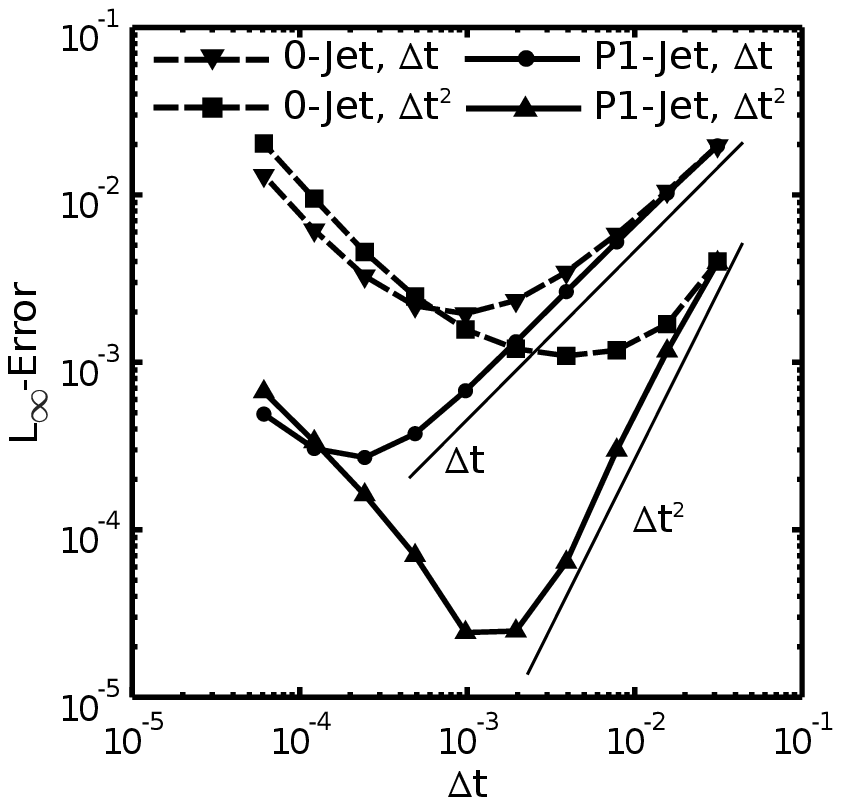}
			\label{fig:comparison_a}					
		}
		\subfigure[$h=0.015625$]{
			\includegraphics[width=0.4\textwidth]{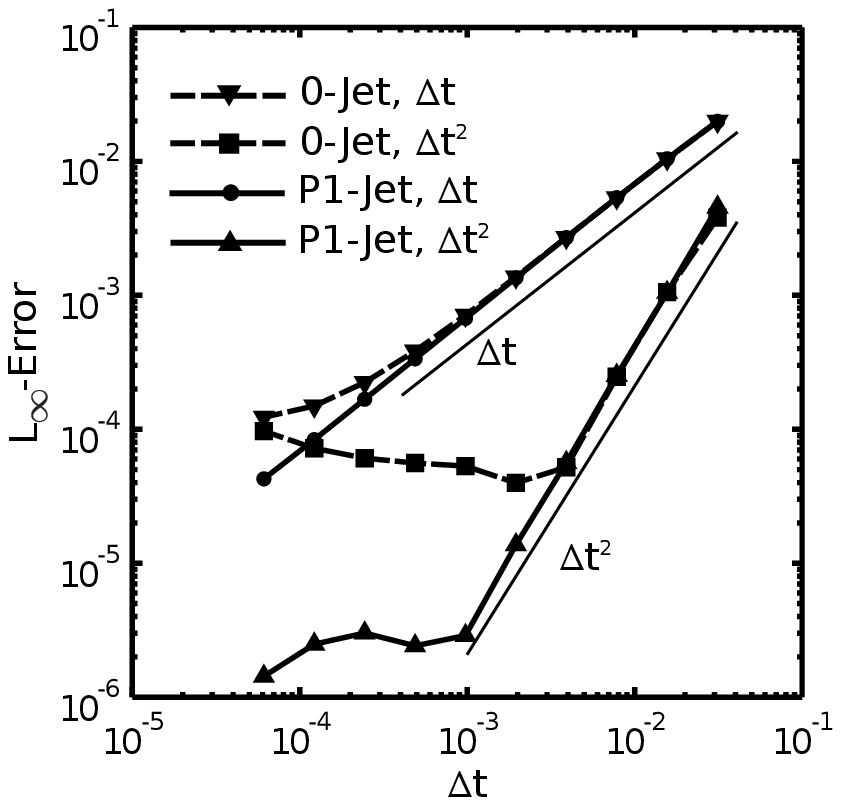}
			\label{fig:comparison_b}				
		}	
		\caption{ Comparison of a 0-Jet with a P1-Jet using both $\mathcal{O}(\Delta t)$ and $\mathcal{O}(\Delta t^2)$ in time methods
			for (a) $h=0.0625$ and (b) $h=0.015625$.}	
		\label{fig:convergence_study}				
	\end{center}
\end{figure}

The first error-convergence results are for fixed grid spacings of $h=0.0625$ and $h=0.015625$, Fig. \ref{fig:convergence_study}.
In this figure both a 0-Jet and a P1-Jet are used with $\mathcal{O}(\Delta t)$ and $\mathcal{O}(\Delta t^2)$ in time methods. Initially
the error converges based on the underlying time integration, $\mathcal{O}(\Delta t^k)$. As the time step decreases the 
$\mathcal{O}(h^2)$ terms begin to dominate and approximately constant error regions are observed, particularly for the 
P1-Jet using second-order time discretization in Fig. \ref{fig:comparison_b}. As the time step further decreases
the dominate term becomes $\mathcal{O}\big(h^{4}/\Delta t\big)$ and the error increases with a further decrease in the time step. 
This is clearly seen in the error plot of the coarser mesh, Fig. \ref{fig:comparison_a}. Both the $\mathcal{O}(\Delta t)$ and $\mathcal{O}(\Delta t^2)$ 
methods have roughly the same error at small time steps as that error is completely dominated by the interpolation error.
For the finer mesh, Fig. \ref{fig:comparison_b}, the integration error is smaller and therefore this increase at small time steps is delayed.
In general the P1-Jet provides more accurate results at moderate and large time steps. Only at small time steps, when the integration error begins to dominate, does the 0-Jet 
have slightly smaller error. The exact cause of this reversal is not known, but could be due to some fortuitous cancelation.

\begin{figure}[ht!]
	\begin{center}		
		\subfigure[0-Jet]{
			\includegraphics[width=0.4\textwidth]{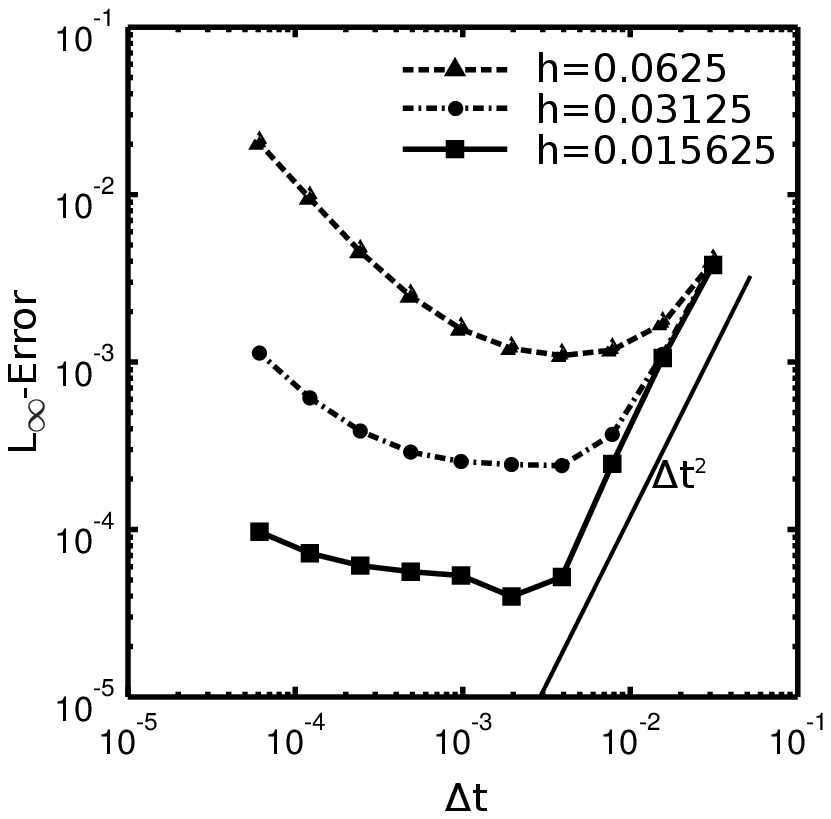}
			\label{fig:timestepping_zero}
		}
		\subfigure[P1-Jet]{
			\includegraphics[width=0.4\textwidth]{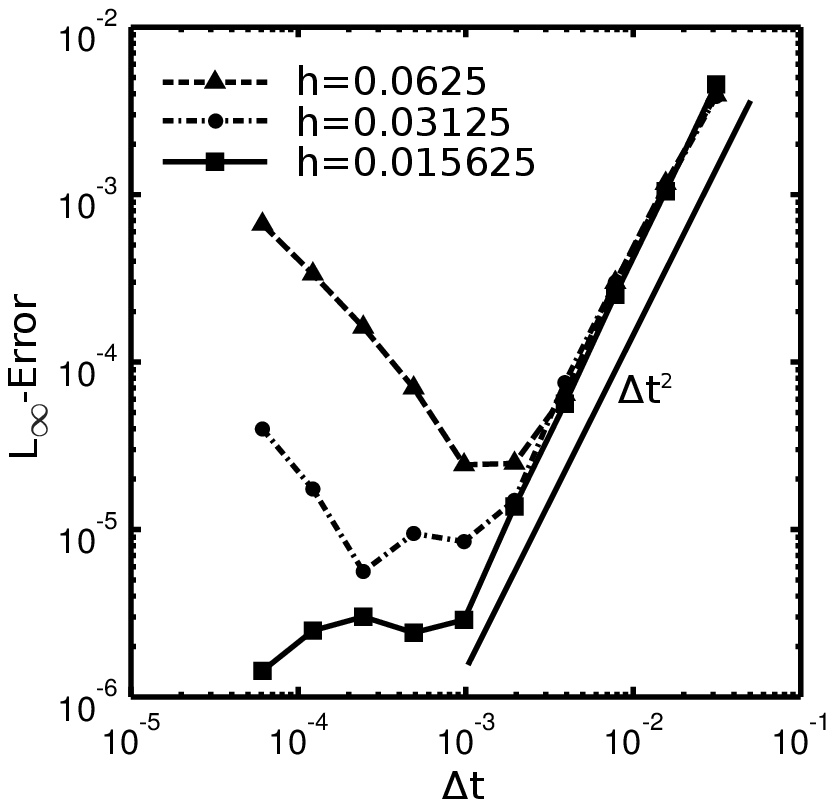}
			\label{fig:timestepping_partial}
		}	
		\caption{ Error dependence on the time step size for different grid sizes. The time-integration scheme is $\mathcal{O}(\Delta t^2)$. }	
		\label{fig:time_stepping}	
	\end{center}
\end{figure}

To further explore the influence of the Jet and of the smoothing operator the error as a function of time step is investigated for different Jets and grid spacings, Fig. \ref{fig:time_stepping}.
At large time steps the error is completely dominated by the time-integration term, and therefore no difference can be observed between the different grid spacings,
particularly for the P1-Jet, Fig. \ref{fig:timestepping_partial}. As the grid spacing is decreased the minimum error over all time steps decreases at a rate of 
$\mathcal{O}(h^2)$, see Table \ref{table:phi_error}.
This indicates that for large time steps any spatial errors are dominated by the smoother operator, not the integration error. This is to be expected so long as $\Delta t>h$,
which results in an overall integration error of $\mathcal{O}(h^3)$.

\begin{table}[ht!]
	\centering		
	\caption{ The smallest $L_\infty$-error of the interface for the circle under mean curvature flow over all time steps. The overall minimum error is dominated
	by the accuracy of the smoothing operation.   }
	\label{table:phi_error}	
	\begin{tabular}{| c | c | c | c | c |}			
		\hline
		$h$ & 0-Jet $L_\infty$ & 0-Jet Order & P1-Jet $L_\infty$ & P1-Jet Order \\
		\hline\hline
		0.0625	& $1.15\times 10^{-3}$ & --- & $2.43\times 10^{-5}$ & ---\\
		0.03125 & $2.45\times 10^{-4}$ & 2.23 & $5.61\times 10^{-6}$ & 2.13\\
		0.015625 & $4.03\times 10^{-5}$ & 2.6 & $1.43\times 10^{-6}$ & 1.96\\
		\hline	
	\end{tabular}		
\end{table}

\subsubsection{Cassini Oval}
\label{subsubsec:3.1.2}

\begin{figure}[ht!]
	\begin{center}		
		  \subfigure[t=0]{
			\includegraphics[width=0.225\textwidth]{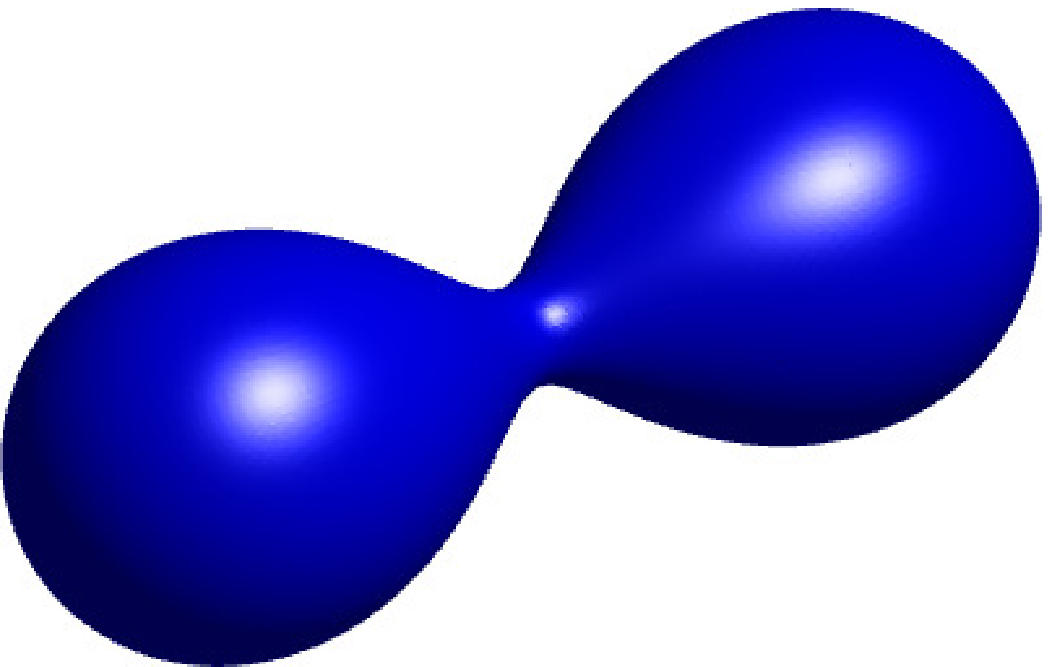}
			\label{fig:3d_cassini_a}
		}
		   \subfigure[t=0.0605]{
			\includegraphics[width=0.225\textwidth]{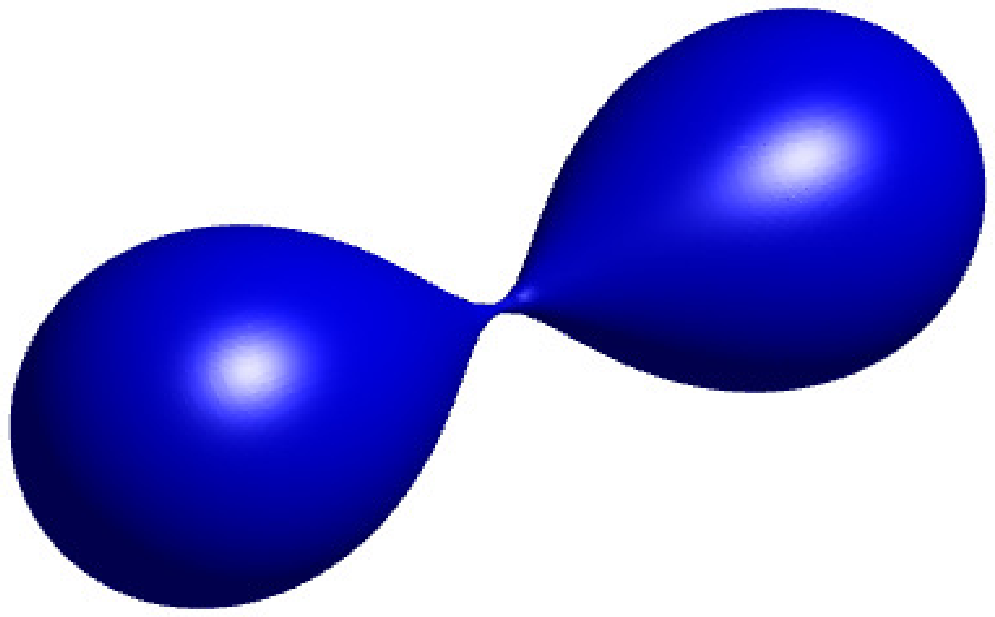}
			\label{fig:3d_cassini_b}
		}		
	       \subfigure[t=0.0645]{
			\includegraphics[width=0.225\textwidth]{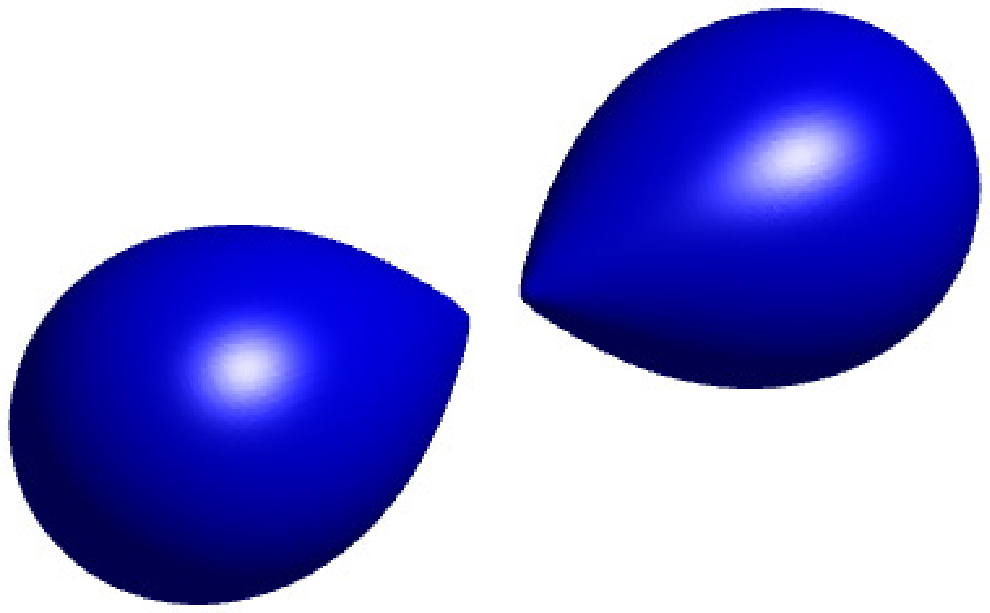}
			\label{fig:3d_cassini_c}
		}
		   \subfigure[t=0.0806]{
			\includegraphics[width=0.225\textwidth]{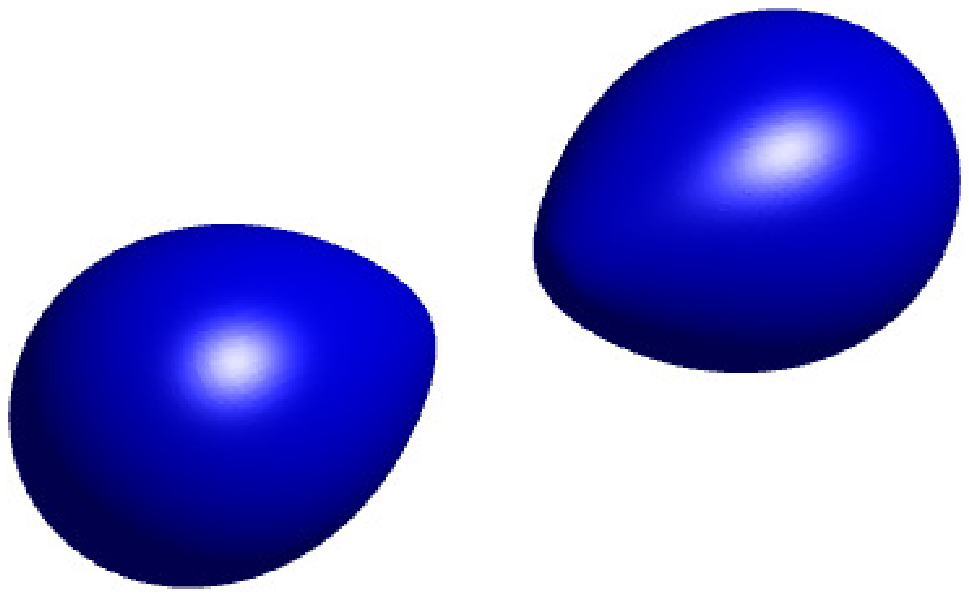}
			\label{fig:3d_cassini_d}
		}			   	
		\caption{Collapse of a three-dimensional Cassini oval under mean curvature flow for $h=0.0634$ and $\Delta t=h^2$ using a P1-Jet with
			second order time discretization.}
		\label{fig:3d_cassini}
	\end{center}				
\end{figure}

The extension of the SemiJet method to three dimensions is demonstrated by considering the collapse
of a Cassini oval under mean-curvature flow, Fig. \ref{fig:3d_cassini}.
The initial shape of the interface is given by $\tilde{\phi}=((x-a)^2+y^2+z^2)((x+a)^2+y^2+z^2)-b^4$
with $a=1.29$ and $b=1.3$ using a grid spacing of $h=0.0634$. This is then replaced by a signed-distance function by reinitialization. 
The surface of the Cassini moves inward and the neck region narrows as the surface shrinks. 
As time goes on the Cassini pinches off at the middle and splits into two separate interfaces. 
The method has been tested using large time steps and works reliably. 
However, the pinch-off dynamics are very rapid and thus to capture the topological changes the time step is chosen to be $\Delta t=h^2$. 
This example also shows that the use of a P1-Jet does not destroy the ability of level set methods to capture topological changes.

\begin{figure}[ht!]
	\begin{center}
		\subfigure[Comparison of a 0-Jet with a P1-Jet]{
			\includegraphics[width=0.45\textwidth]{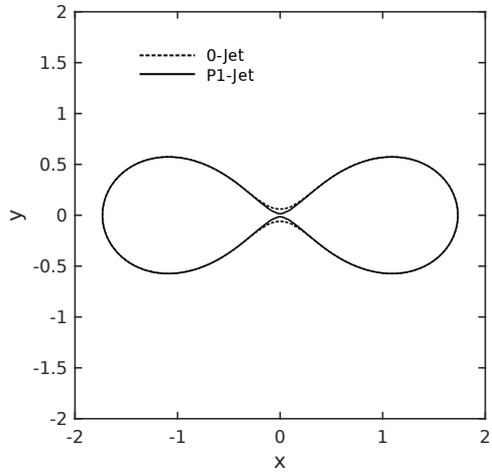}
			\label{fig:crosssection}				
		} \qquad
		\subfigure[Close-up view of the neck region. 0-Jet under-reports the curvature without the higher accuracy afforded by P1-Jet]{
			\includegraphics[width=0.45\textwidth]{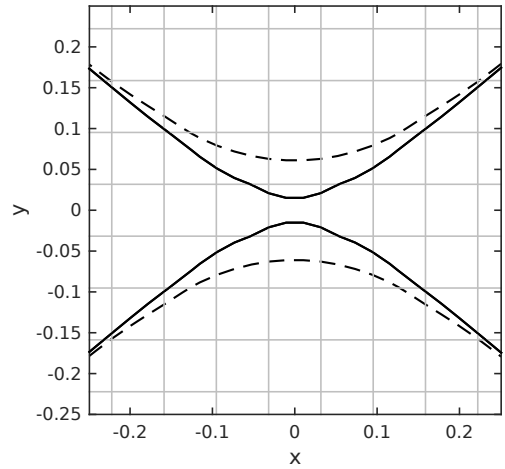}
			\label{fig:zoom}				
		}
	\end{center}	
	\caption{Cross-section of the $x-y$-plane of the Cassini Oval using a 0-Jet and P1-Jet at a time of $t=0.0605$. The neck region of the 0-Jet is 
		thicker than the P1-Jet.}
	\label{fig:2d_cassini}				
\end{figure}

A comparison between the 0-Jet and P1-Jet for a collapsing Cassini oval is shown in Fig. \ref{fig:2d_cassini}. As can be seen
the neck using the 0-Jet is thicker than the P1-Jet at the same time. This indicates that
the 0-Jet is under-reporting the curvature, leading to a delay in splitting time. 

\begin{figure}[ht!]
	\begin{center}
		\includegraphics[width=3in]{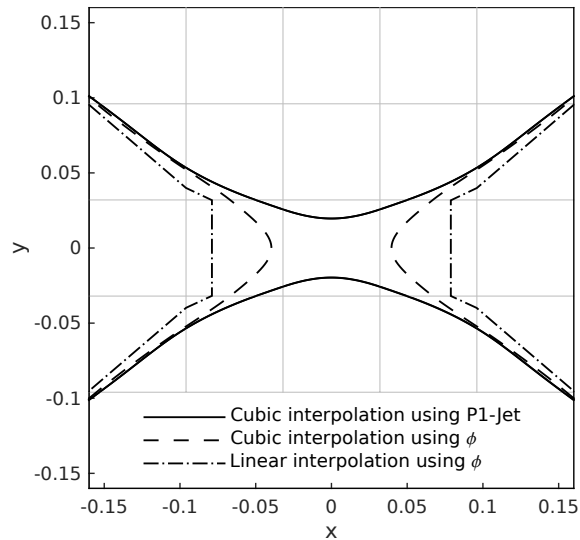}
		\caption{Close-up view of the Cassini Oval at t = 0.0605 showing sub-grid structure. The neck region falls within a grid cell (represented by solid black lines) which is 
		captured by P1-Jet through the locally available gradient information}
		\label{fig:P1jet_compare}
	\end{center}
\end{figure}

In Fig. \ref{fig:P1jet_compare} the ability of the P1-Jet to capture sub-grid information is demonstrated. Here, a comparison is made 
between interfaces constructed with three different interpolation approaches for the collapsing Cassini oval at a time of $t=0.0605$.
The solid black lines represents the interface constructed using the P1-Jet information to determine the cubic interpolant. 
This result is the same as Fig. \ref{fig:zoom}. This is denoted as a sub-grid structure due to the fact that two interface crossings occur
on a cell edge where $\phi>0$ on both of the cell corners. Compare this to a cubic interpolant determined using only
the level set values and approximating the derivatives from finite difference approximations, the dashed line in Fig. \ref{fig:P1jet_compare}.
In this case the center region is completely split. Finally, consider the case of using linear interpolation with only the level set values,
the dash-dot line in Fig. \ref{fig:P1jet_compare}. The description of the interface using linear interpolation is even worse, with only
straight lines being described. This demonstrates that the tracking of level set derivatives is beneficial in describing the interface.

\subsection{Volume-Conserving Mean Curvature Flow}
\label{subsec:3.2}
In volume-conserving mean curvature flow the velocity of the interface is given by $\vec{u}=-(\kappa-\kappa_{avg})\vec{n}$ where $\kappa$ is the mean 
curvature, $\vec{n}$ is the outward normal to the interface and $\kappa_{avg}$ is the average mean curvature on the interface given by
\begin{lequation}
	\kappa_{avg}=\frac{\int_{\Gamma} \kappa\;\mathrm{d}a}{\int_{\Gamma}\;\mathrm{d}a}\approx\frac{\int_\Omega \kappa\;\delta(\phi)\;\|\nabla\phi\|\;\mathrm{d}\vec x}{\int_\Omega \delta(\phi)\;\|\nabla\phi\|\; \mathrm{d}\vec x},
\end{lequation}
where $\Omega$ is the computational domain and $\delta(\phi)$ is a Dirac delta function \cite{Towers2008}.
Under this velocity field the interface will evolve into a single circle in 2D and sphere in 3D, where the average curvature equals the local curvature.

\begin{figure}[ht!]
	\begin{center}
		   \subfigure[t=0]{
			\includegraphics[width=0.225\textwidth]{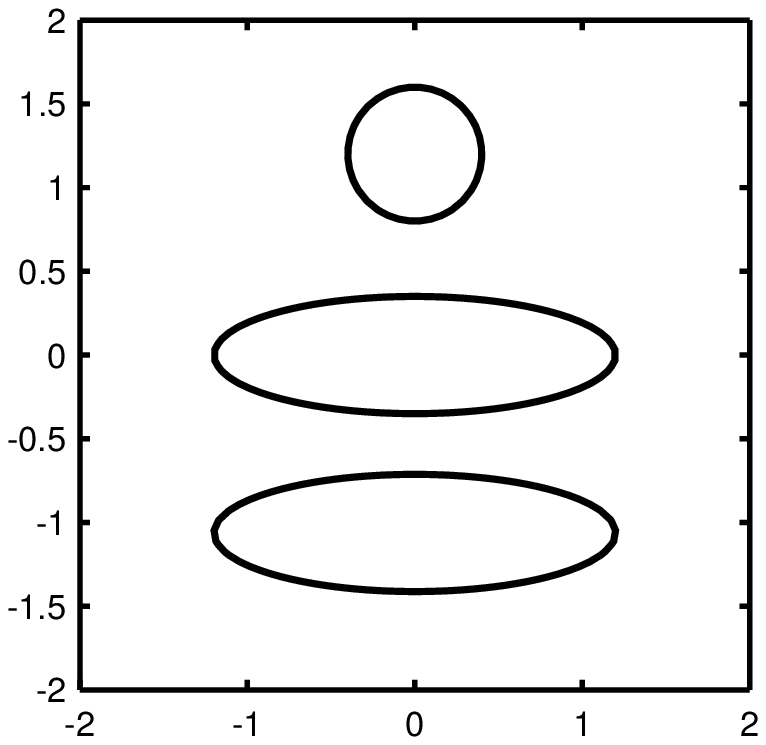}
			\label{fig:area_conserving_zero_a}
		}
		   \subfigure[t=0.10]{
			\includegraphics[width=0.225\textwidth]{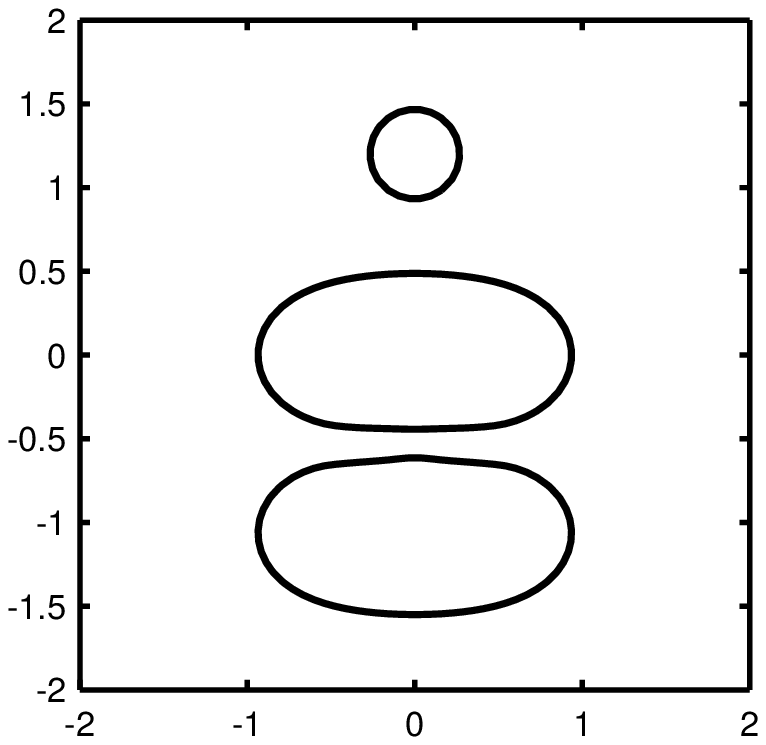}
			\label{fig:area_conserving_zero_b}
		}
		  \subfigure[t=0.12]{
			\includegraphics[width=0.225\textwidth]{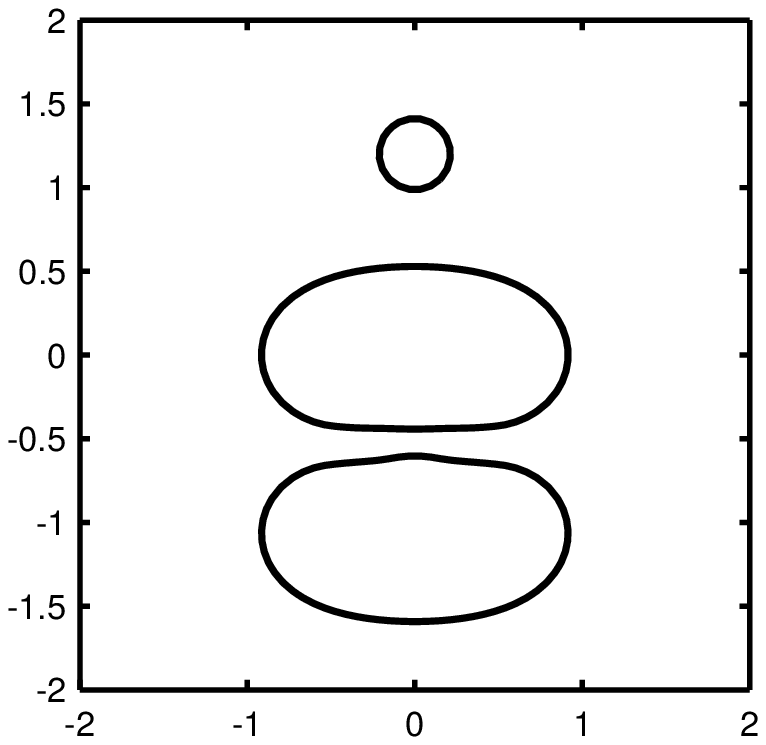}
			\label{fig:area_conserving_zero_c}
		} 
		   \subfigure[t=0.15]{
			\includegraphics[width=0.225\textwidth]{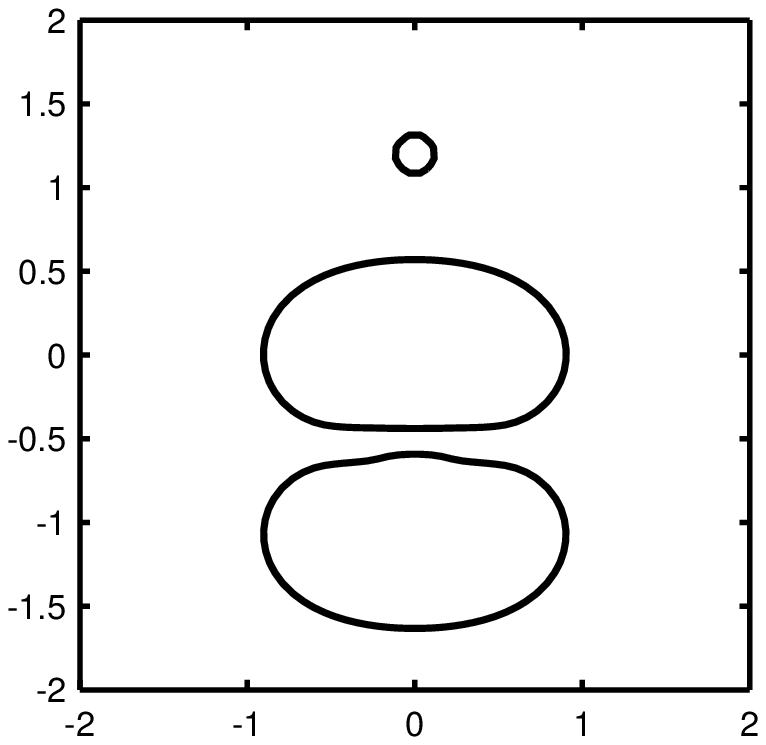}
			\label{fig:area_conserving_zero_d}
		} 		
		\caption{ Volume-conserving mean curvature flow in two dimensions using a 0-Jet. The grid spacing is $h=0.0634$ while the time step is $\Delta t=0.5h^2$.
			At $t=0.15$ the volume change is 1.35\%.
		}	
		\label{fig:area_conserving_zero_jet}
	\end{center}
\end{figure}

\begin{figure}[ht!]
	\begin{center}
		   \subfigure[t=0]{
			\includegraphics[width=0.225\textwidth]{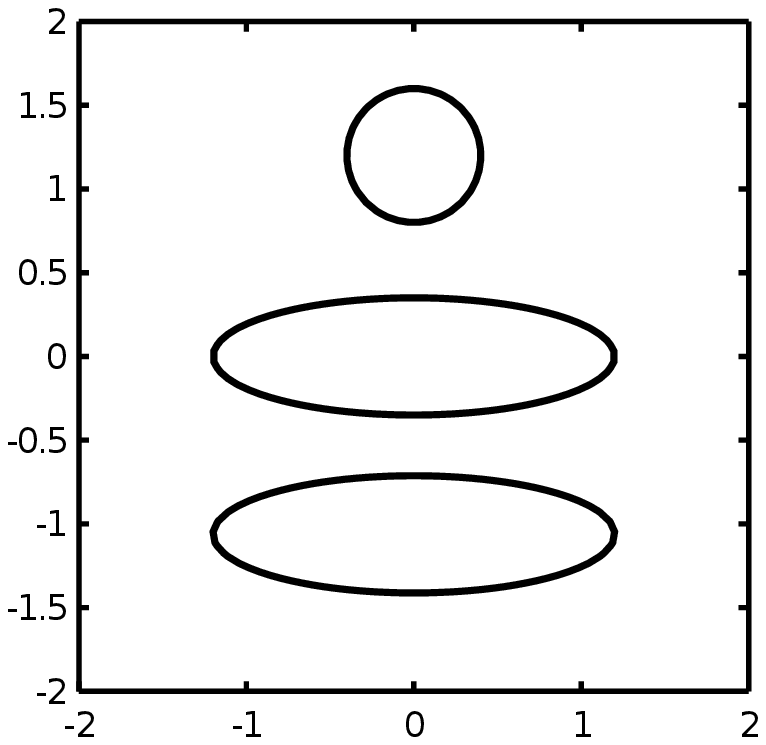}
			\label{fig:area_conserving_a}
		}
		   \subfigure[t=0.10]{
			\includegraphics[width=0.225\textwidth]{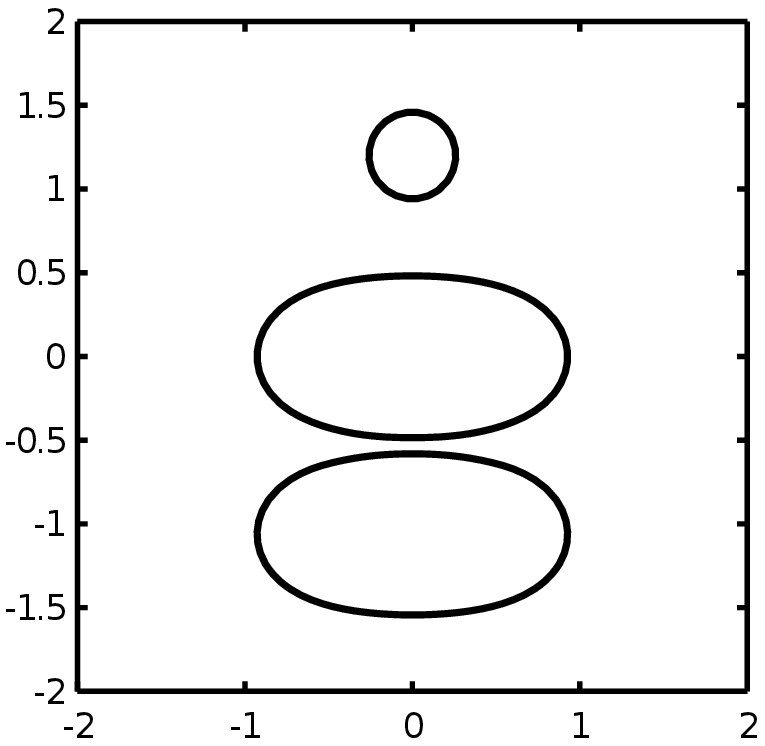}
			\label{fig:area_conserving_b}
		}
		  \subfigure[t=0.12]{
			\includegraphics[width=0.225\textwidth]{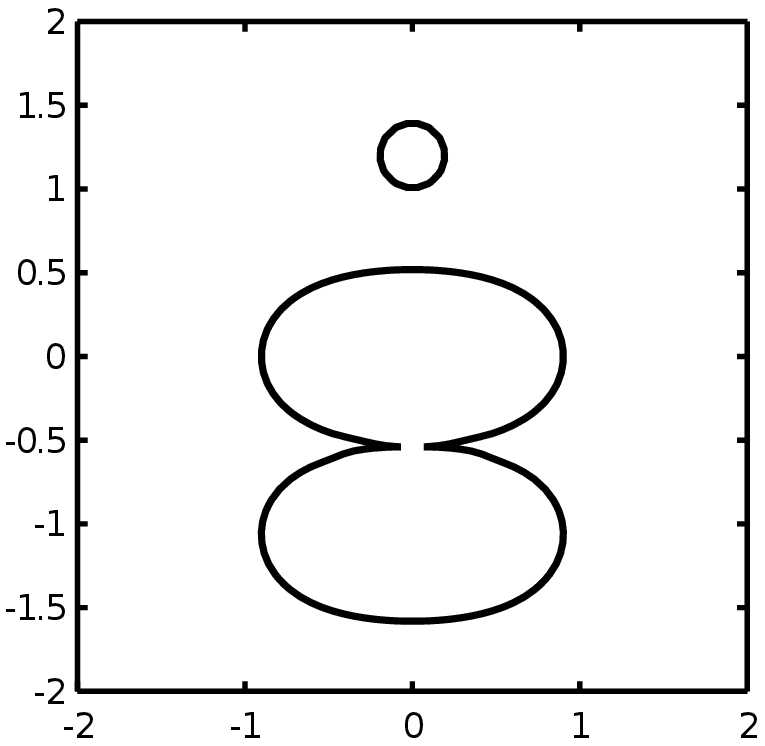}
			\label{fig:area_conserving_c}
		} 
		   \subfigure[t=0.15]{
			\includegraphics[width=0.225\textwidth]{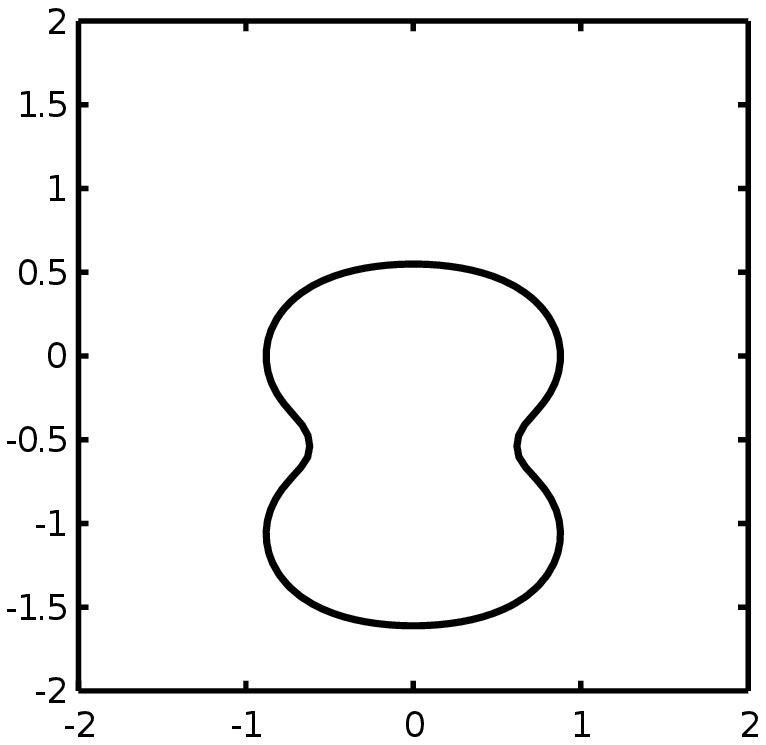}
			\label{fig:area_conserving_d}
		} 		
		\caption{ Volume-conserving mean curvature flow in two dimensions using a P1-Jet. The grid spacing is $h=0.0634$ while the time step is $\Delta t=0.5h^2$. 
		At $t=0.15$ the volume change is 1.28\%.
		}	
		\label{fig:area_conserving}
	\end{center}
\end{figure}

The first test case is a two-dimensional system containing two ellipses and a circle subjected to volume conserving mean curvature flow.
As the curvature is a quantity local to each interface, ideally all interfaces would evolve independently until they merge. 
Provided enough time the final shape of the interface will become a circle. Results for both a 0-Jet, Fig. \ref{fig:area_conserving_zero_jet},
and P1-Jet, Fig. \ref{fig:area_conserving}, are provided.
The volume loss for the 0-Jet is 1.35\% while that for the P1-Jet is 1.28\%. 

While the volume conservation is similar it is clear that the 0-Jet 
result is incorrect.
It has been previously observed that the non-local nature of the semi-implicit level set method results in the incorrect behavior of near-by interfaces 
\cite{salac2008local, Smereka2003}. 
This is due to the difficulty in calculating the curvature when two interfaces are near each other \cite{macklin2006improved, lervaag2013calculation}.
This can be clearly observed in the 0-Jet time-evolution, where the interfaces are clearly distorted and have not yet merged at $t=0.15$, Fig. \ref{fig:area_conserving_zero_jet}.
Contrast that with the P1-Jet, where much less interface distortion is observed and the interfaces are fully merged at $t=0.15$. 
While the curvature calculation is non-local to the interface using the curvature calculation described earlier, the higher accuracy, as demonstrated by 
the results in previous sections, reduces the influence of nearby interfaces.

\begin{figure}[ht!]
	\begin{center}
		
			\subfigure[$t=0$]{				     
				\includegraphics[width=0.225\textwidth]{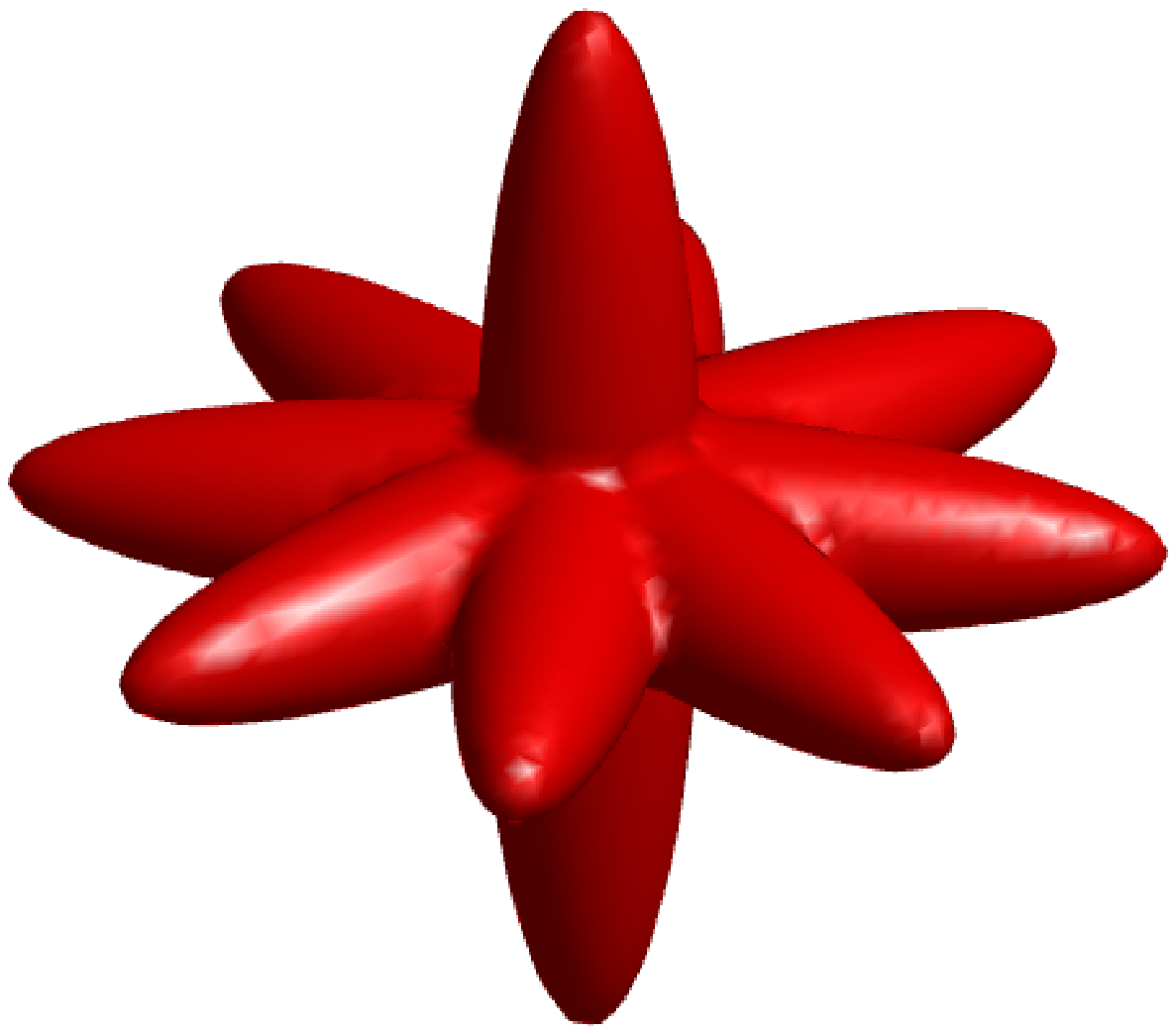}
			}
			\subfigure[$t=0.05$]{				     
				\includegraphics[width=0.225\textwidth]{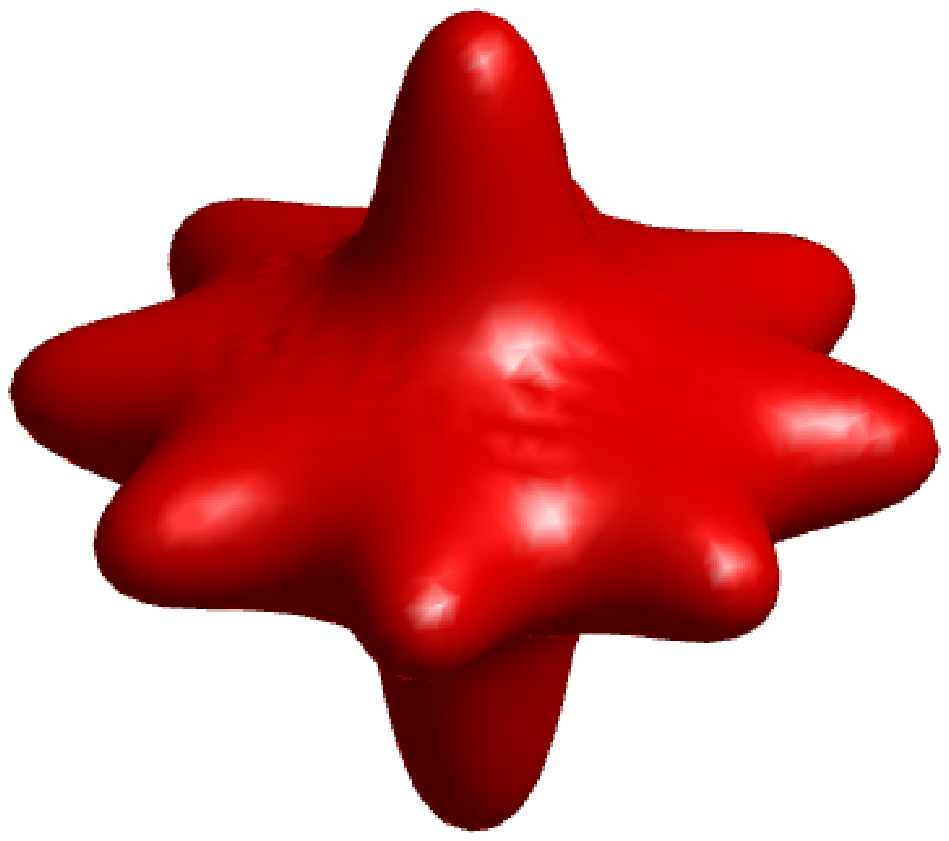}
			}
			\subfigure[$t=0.11$]{				     
				\includegraphics[width=0.225\textwidth]{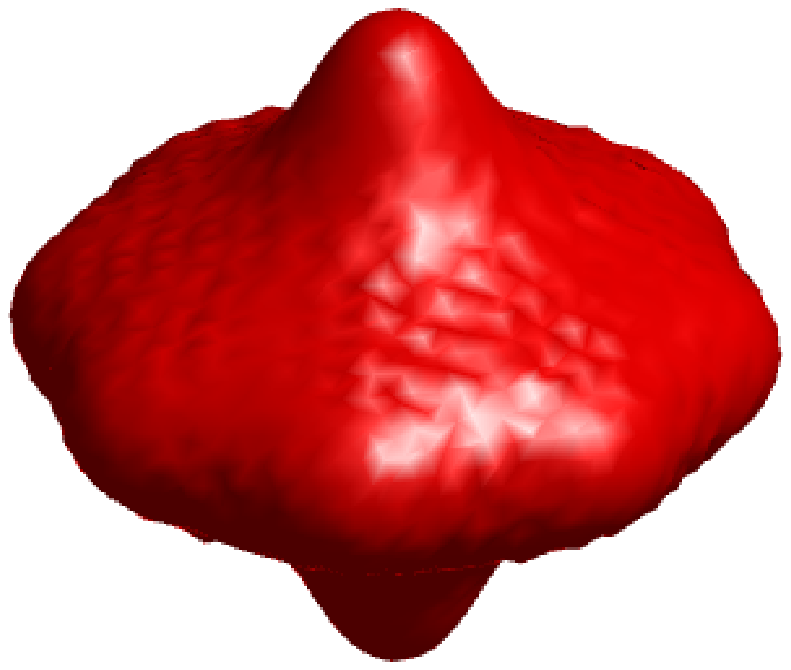}
			}
			\subfigure[$t=0.58$]{				     
				\includegraphics[width=0.225\textwidth]{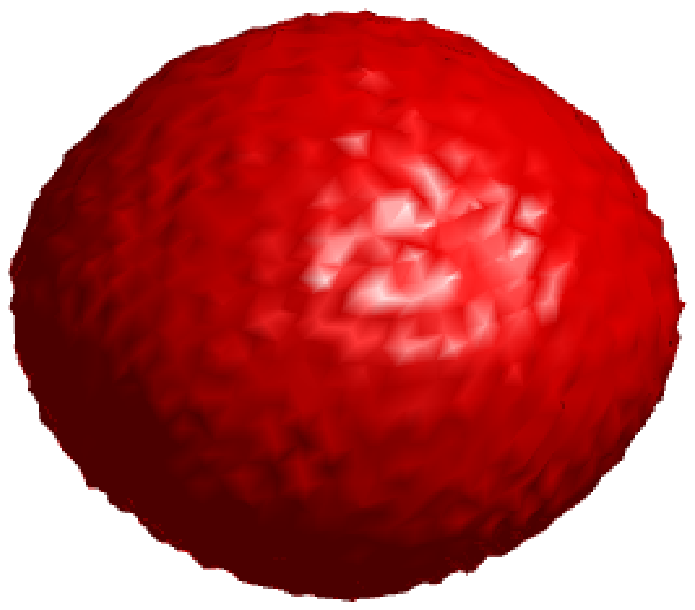}
			}\\		
			\subfigure[$t=0$]{				     
				\includegraphics[width=0.225\textwidth]{star_partial1}
			}
			\subfigure[$t=0.05$]{				     
				\includegraphics[width=0.225\textwidth]{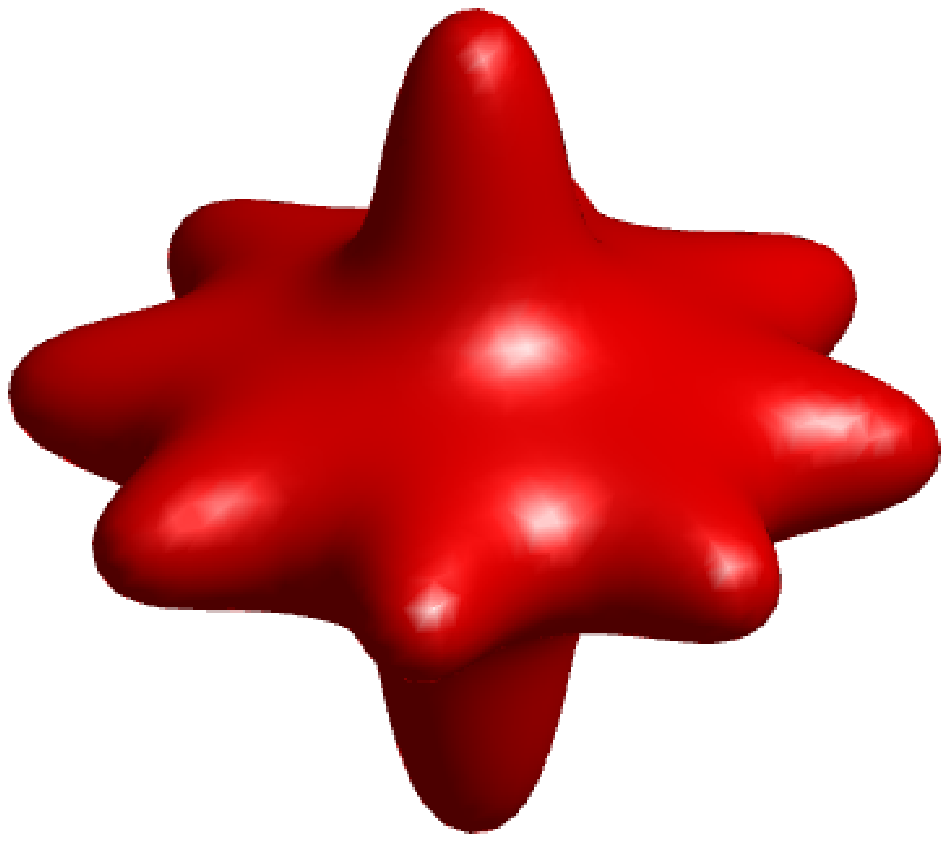}
			}
			\subfigure[$t=0.11$]{				     
				\includegraphics[width=0.225\textwidth]{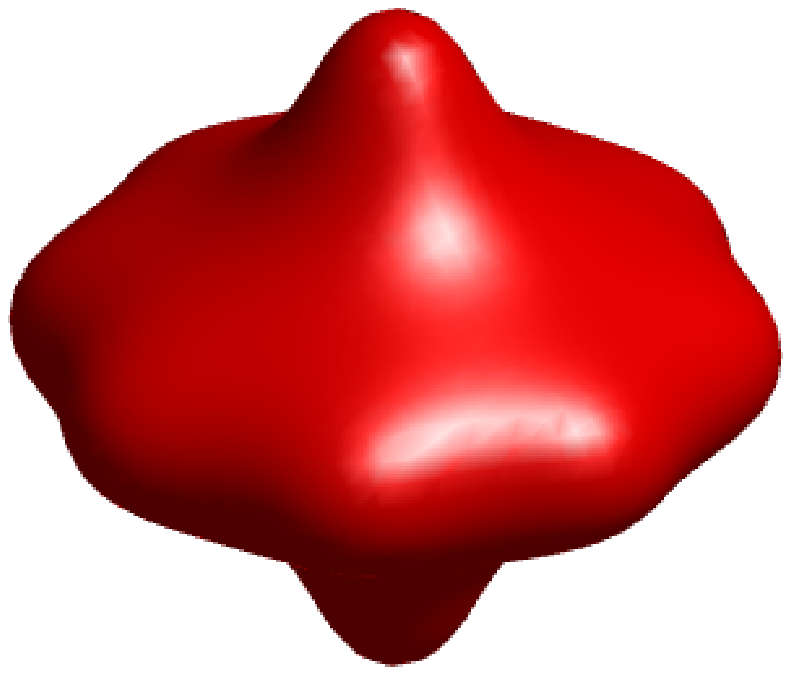}
			}
			\subfigure[$t=0.58$]{				     
				\includegraphics[width=0.225\textwidth]{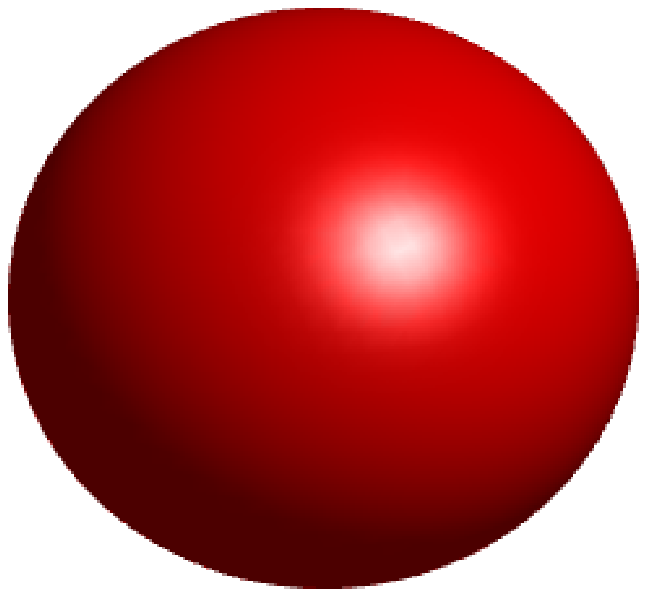}
			}								
		\caption{A comparison of the original Jet ($\beta=0$) method, (a)-(d), with the SemiJet, (e)-(h), showing the volume-conserving
		mean curvature flow of a three dimensional star	with $\Delta t=1.5h^2$ and $h=0.0625$. Both solutions use a P1-Jet.
		Using the original Jet method oscillations on the interface grow over time, leading 
		to the eventual failure of the simulation. The SemiJet scheme, on the other hand, maintains a smooth solution at all times. 
		Using the SemiJet volume loss is about 0.37\% at $t=0.58$.}
		\label{fig:3d_star_comparison}
	\end{center}
\end{figure}

A comparison between the original, explicit Jet scheme is also made to the presented SemiJet method by considering the volume-conserving 
mean-curvature flow of a three-dimensional star shape, Fig. \ref{fig:3d_star_comparison}.
Under this motion any given initial shape will eventually become a sphere. 
Using the original Jet scheme oscillations begin to appear on the interface, which will lead to numerical failure due to incorrect curvature calculation.
Numerical experiments show that the instability of the original Jet scheme becomes more severe over time and the solution eventually breaks down.
The presented SemiJet method is able to maintain a smooth interface through the entire time with the same time step with only a volume change of 0.37\%.

\subsection{Qualitative Stability Analysis}
\label{subsec:stability}
In this section the stability of the SemiJet method (P1-Jet with $\beta=0.5$) 
is compared to the original semi-implicit level set method (0-Jet with $\beta=0.5$) and the original Jet scheme (P1-Jet with $\beta=0$).
Consider a two-dimensional circle under volume-conserving mean curvature flow, $\vec{u}=-(\kappa-\kappa_{avg})\vec{n}$. In the ideal case the interface will not move due to 
the local curvature being equal to the average curvature over the interface. This situation can be used to investigate the stability of
the underlying numerical scheme. As the interface should not move any motion of the interface will be due to errors
when calculating the curvature and advecting the interface. If the time step is too large for a given numerical method these errors will compound resulting
in the eventual failure of the simulation.

To test the stability various numerical simulations are carried out to determine the maximum stable time step for various combinations of Jet-order, time step
accuracy order, and amount of smoothing. While the amount of volume change during the course of a simulation is small over reasonable time steps, as demonstrated above, during large 
time steps the volume change can become excessive. To minimize the influence of this volume change on the determination of the maximum stable time step a volume correction procedure is 
performed. After every time step the level set is shifted by $\left(V_n-V_0\right)/A_n$, where $V_n$ is the current volume, $V_0$ is the initial volume, and $A_n$ is the current 
area of the interface. In 2D these are replaced with the current area, initial area, and current interface length, respectively.

\begin{figure}[ht!]
	\begin{center}
		\includegraphics[width=3in]{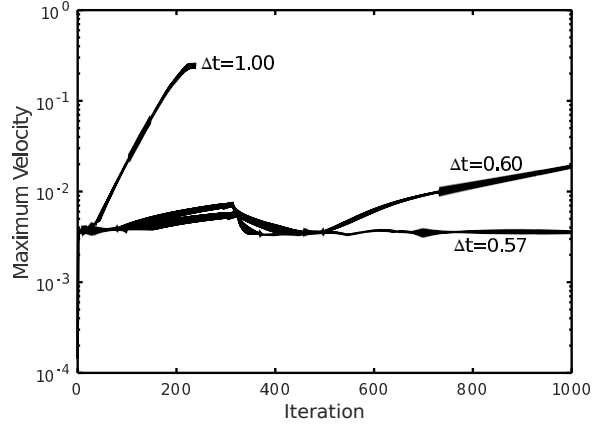}
		\caption{The maximum velocity using three different time steps using a 0-Jet with $\beta=0.5$ for a unit circle under volume-conserving
		mean curvature flow.
		Three types of behaviors are shown. An unstable time step ($\Delta t=1.00$) does 
		not complete 1000 iterations. A semi-stable time step ($\Delta t=0.60$) completes the 1000 iterations but the velocity is slowly increasing over time. A stable 
		time step ($\Delta t=0.57$) completes 1000 iterations and the velocity remains constant.}
		\label{fig:MaxVelocity}
	\end{center}
\end{figure}

A time step is determined to result in a stable simulation if the following conditions are satisfied: 1) the simulation
can run for 1000 time steps and 2) the difference between the maximum and minimum velocity over the last 100 iterations is less than 10\% of the minimum velocity over 
the 1000 total time steps. Sample velocity profiles over 1000 iterations for a 0-Jet using $\beta=0.5$ with second-order accurate time discretization for three 
different time steps is shown in Fig. \ref{fig:MaxVelocity}. There are three types of behaviors observed. The first is an unstable time step, where the simulation
is not able to complete 1000 iterations before failing due to excessive numerical errors. The second type is a semi-stable time step, where the simulation
completes 1000 iterations but the velocity is slowly increasing. This indicates that the simulation will fail at some point in the future. Finally, a stable time 
is one which completes all 1000 iterations and the velocity is stable. A stable time step will not have zero velocity; instead any small movement of the interface
is adjusted by the volume-correction procedure. 

The maximum possible stable time step is reported in Table \ref{table:maxDT} for a two-dimensional unit circle with volume-conserving mean curvature flow 
on a grid with $h=0.0625$.
While this determination of a stable time step is qualitative in nature it does provide insight into the overall
stability properties of the SemiJet method.
It is clear that the explicit scheme, given by $\beta=0$, is restricted by the standard CFL condition of $\Delta t\leq h^2$.
The application of a smoothing operator, given by $\beta=0.5$, relaxes this constraint.
In the case of a P1-Jet and first-order time discretization the maximum possible
time step is over 400-times that of $h^2$. 

\begin{table}[ht!]
	\centering		
	\caption{ The largest stable time step for various combinations of the Jet, time step order, and smoothing on a grid with $h=0.0625$. The explicit Jet method is given by $\beta=0$ while
	the SemiJet method is given by $\beta=0.5$.}
	\label{table:maxDT}	
	\begin{tabular}{|c|c|c|c|c|c|}
		\cline{3-6} 
		\multicolumn{1}{c}{} &  & \multicolumn{2}{c|}{$\beta=0$} & \multicolumn{2}{c|}{$\beta=0.5$}\\
		\hline 
		Jet Order & Time Order & $\Delta t_{max}$ & $\Delta t_{max}/h^2$ & $\Delta t_{max}$ & $\Delta t_{max}/h^2$ \\
		\hline\hline	
		0 & 1 &  0.00205 & 0.52 &  0.98 & 251\\ 	
		P1 & 1 &  0.00375 & 0.96 &  1.71 & 438\\ 	
		0 & 2 &  0.00135 & 0.35 &  0.57 & 146\\ 	
		P1 & 2 &  0.00245 & 0.63 &  0.79 & 179\\ 	
		\hline	
	\end{tabular}		
\end{table}

It is interesting to note that in all cases the use of a P1-Jet results in a more stable scheme than the 0-Jet. This holds true for both the explicit case, $\beta=0$, and the 
semi-implicit case, $\beta=0.5$. It is suspected that the averaging of the derivatives when calculating the curvature in the P1-Jet case is partially responsible for this increase in stability.
It should also be noted that in general a first-order time discretization scheme is more stable than the second-order schemes. This is due to the 
extrapolation step during the second-order time discretization. Any oscillations in the velocity field will be magnified when extrapolating into the future.

\subsection{Computational Cost Comparison}

In this section the computational cost of the proposed SemiJet method is compared to a standard fifth order upwind WENO method. 
The \newtext{first} test case considered here is the collapse of a two-dimensional unit circle under mean curvature flow, $\vec{u}=-\kappa\vec{n}$. The system is allowed to run until $t=0.375$, 
at which time the error is computed, see Sec. \ref{subsec:mean-curvature-flow} for more details. 
\newtext{The domain spans the $[-2,2]^2$ square.}
During the simulation the wall-clock time spent in the evolution of the level-set or Jet is determined. Note that 
operations common to different cases, such as reinitialization, are not included in the timings. This was done to ensure that an accurate comparison
between methods is performed. All simulations were performed using 8 AMD Opteron 6320 \newtext{cores}.

The WENO method is discretized using a second-order backward-finite-difference method in time and second-order extrapolation of the velocity in time:
\begin{lequation}
	\dfrac{3 \phi_{n+1}-2\phi_n+\phi_{n-1}}{2\Delta t}+2\vec{u}_{n}\cdot\nabla\phi_{n}-\vec{u}_{n-1}\cdot\nabla\phi_{n-1}=0,
\end{lequation}
where the gradients are discretized using a fifth-order upwind WENO scheme \cite{Jiang2000}. Note that when using the WENO scheme,
only the level-set field is tracked. The Jet (with $\beta=0$) and SemiJet (with $\beta=0.5$) both 
utilize a P1-Jet. As $\beta=0$ for the Jet, it is an explicit scheme in time, and thus both the WENO and Jet schemes will be under 
strict stability constraints, see Sec. \ref{subsec:stability} for a qualitative stability analysis of of the Jet scheme. It is expected
that WENO will obey similar restrictions.
The results are presented in Table \ref{table:2D_Timings}. For the WENO and Jet methods the time step is set to $0.25h^2$, while with the SemiJet
method a larger time step can be used.

\newtext{
\begin{table}[ht!]
	\centering		
	\caption{ The $L_\infty$ error and wall-clock time required to model a two-dimensional circle collapsing under mean-curvature flow until a time of $t=0.375$.\vspace{2ex}}
	\label{table:2D_Timings}	
	\begin{tabular}{|c|c|c|c|c|c|}		
		\hline 
		Method & Grid Size & Time Steps & $\Delta t/h^2$ & $L_\infty$ Error & Elapsed Time (s) \\
		\hline\hline	
		WENO & 257 & 6144 & 0.25 		& $6.81\times 10^{-3}$ 	& 99.53\\
		WENO & 129 & 1536 & 0.25 		& $1.33\times 10^{-2}$ 	& 6.38 \\
		WENO & 65 & 384 & 0.25 			& $2.54\times 10^{-2}$ 	& 0.42 \\
        \hline
       Jet & 257 & 6144 & 0.25 			& $2.03\times 10^{-6}$ 	& 1527 \\
       Jet & 129 & 1536 & 0.25 			& $1.06\times 10^{-5}$ 	& 91.95 \\
       Jet & 65 & 384 & 0.25 			& $4.68\times 10^{-5}$ 	& 5.95 \\
        \hline     
		SemiJet & 257 & 6144 & 0.25 	& $2.41\times 10^{-6}$ 	& 2377  \\
		SemiJet & 257 & 250 & 6.14 		& $6.82\times 10^{-6}$ 	& 98.36 \\
		SemiJet & 257 & 24 & 64 		& $1.04\times 10^{-3}$ 	& 10.50 \\
		\hline
		SemiJet & 129 & 1536 & 0.25 	& $1.94\times 10^{-5}$ 	& 148.8 \\
		SemiJet & 129 & 48 & 8 			& $2.91\times 10^{-4}$ 	& 4.76 \\
		SemiJet & 129 & 24 & 16 		& $1.07\times 10^{-3}$ 	& 2.45 \\
		\hline
		SemiJet & 65 & 384 & 0.25 		& $6.87\times 10^{-5}$ 	& 9.91 \\
		SemiJet & 65 & 12 & 8 			& $3.91\times 10^{-3}$ 	& 0.35 \\
		\hline
		SemiJet & 33 & 12 & 2 			& $4.33\times 10^{-3}$ 	& 0.12 \\		
		\hline	
	\end{tabular}		
\end{table}
}
It is clear from the results that the Jet and SemiJet are much more costly than the WENO 
scheme for the same grid size and time step, with the Jet being approximately \newtext{15} times more
expensive and the SemiJet approximately \newtext{23} times more expensive. This is due to the fact that 
for every grid point, four sub-grid points (in 2D) must be updated. Each of these updates 
requires not only determining the departure location, but also the evaluation of multiple
cubic interpolants. Note that no effort was made to optimize this interpolation, and 
thus some time savings could be gained there.

The results also indicate, though, that the quality of the solution, as measured by the 
$L_\infty$ error, is much better for the Jet and SemiJet schemes as compared to the 
WENO method. For the same grid size and time step, the Jet and SemiJet schemes
are between two- and three-orders of magnitude more accurate than the WENO scheme. 
For the SemiJet scheme, this allows for an increase in the time step to still
produce an accurate result. For example, using a grid \newtext{size of $257^2$ (grid spacing of $1.5625\times 10^{-2}$)} the WENO result 
has an error of $6.81\times 10^{-3}$. The same computational cost ($\sim 100$ s)
can be achieved with the SemiJet scheme using a time step 24 times larger than the WENO case.
Even with this larger time step, the error for the SemiJet case is three orders of magnitude smaller, $6.82\times 10^{-6}$ for the SemiJet scheme. 
Using the SemiJet scheme with a time step
256 times larger ($\Delta t=1.563\times 10^{-2}$) than the WENO one for the same grid results in a scheme which is 6.5 times more accurate
but only takes one-tenth of the time, 10.5 s versus 99.53 s.

This behavior, that at the same computational cost the SemiJet scheme is 
much more accurate than the WENO scheme is consistent across multiple grid sizes. 
In fact, using an extremely coarse grid spacing of 0.125, \newtext{corresponding to a $33^2$ grid}, with a very large
time step still results in a more accurate result than the WENO 
scheme with a fine grid and very small time step, in a computational
time which is 830 times faster. This behavior is consistent with that shown by
Chidyagwai et al. for linear advection \cite{Chidyagwai}.

\newtext{
The second test case is mean curvature flow in three dimensions, $\vec{u}=-0.5\kappa\vec{n}$, where $\kappa$ represents the sum
of the principle curvatures. In this case the size of the grid
is fixed at $129^3$ over a cube spanning $[-2,2]^3$. These simulations were performed using 128 AMD Opteron 6320 cores, with the 
results seen in Table \ref{table:3D_Timings}. 
}
\begin{table}[ht!]
	\centering		
	\caption{ The $L_\infty$ error and wall-clock time required to model a three-dimensional sphere collapsing under mean-curvature flow until a time of $t=0.375$.
			All cases use a $129^3$ grid. \vspace{2ex}}
	\label{table:3D_Timings}	
	\begin{tabular}{|c|c|c|c|c|}		
		\hline 
		Method & Time Steps & $\Delta t/h^2$ & $L_\infty$ Error & Elapsed Time (s)  \\
		\hline\hline	
		WENO 	& 1536 	& 0.25 		& $2.43\times 10^{-2}$ 	& 77.33 \\
		\hline
		Jet 	& 1536 	& 0.25 		& $1.14\times 10^{-4}$ 	& 2022 \\
		\hline
		SemiJet & 1536 	& 0.25 		& $1.21\times 10^{-4}$ 	& 3417 \\
		SemiJet & 48 	& 8 		& $2.73\times 10^{-4}$ 	& 116.9 \\
		SemiJet & 30 	& 12.8 		& $6.29\times 10^{-4}$ 	& 73 \\
		SemiJet & 24 	& 16 		& $9.65\times 10^{-4}$ 	& 60.15 \\
		SemiJet & 12 	& 32 		& $3.63\times 10^{-3}$ 	& 32.44 \\
		\hline	
	\end{tabular}		
\end{table}
\newtext{
The results for the three-dimensional case are similar to those for the two-dimensional one. Here, the Jet is approximately 26 times
more expensive and the SemiJet is approximately 44 times more expensive than the WENO scheme for the same grid size. In three-dimensions
there are eight sub-grid points that need to be updated, which explains the increase over the two-dimensional results. 
Similar to the two-dimensional case, the
Jet and SemiJet schemes are two orders of magnitude more accurate than the WENO scheme at the same grid size and time-step. 
Using the SemiJet method with 30 time steps, corresponding
to $\Delta t=0.0125$, results in the same computational cost as the WENO method but a result which is 40 times more accurate. In fact, using an extremely coarse $33^3$ mesh
with a time step of $\Delta t=0.125$ on two cores results in an elapsed time of 21 s and an error of $3.05\times 10^{-3}$, which is both faster 
and more accurate than the WENO case, using much fewer computational resources.
}

\section{Summary and Conclusion}
\label{sec:conclusion}
In this work a robust numerical method, the SemiJet, is presented for moving interface problems with numerically-stiff
velocity fields. The method builds upon the original Jet 
scheme developed by Seibold et al. and the semi-implicit level set method of Smereka. The influence of the semi-implicit
smoothing operation is captured by defining a smoothing source term, which is then used when updating the derivative fields in the level set Jet. 
The results demonstrate that the use of a Jet increases both the accuracy and stability of the simulations over both the original Jet and original
semi-implicit level set methods. Unlike the previous SIGALS method only one smoothing operation is required 
and the use of the $\epsilon$-finite different method allows 
for straight-forward extension to higher order jet schemes. While the SemiJet method is computationally more expensive than a WENO scheme
for the same grid and time step,
the added stability and accuracy allows for smaller grid sizes and larger time steps to be used, which results in methods with overall 
higher efficiency and accuracy.

\section*{Acknowledgments}
	This work has been supported by the National Science Foundation through the Division of Chemical, Bioengineering, Environmental, and Transport Systems Grant \#1253739.
	The authors would also like to thank Dr. Paul Bauman for discussions on this work.

\section*{References}
\bibliography{semijet}
\end{document}